\documentclass[11pt]{article}
\usepackage{bbm,theorem,latexsym,graphicx}

\let\phi\varphi
\let\theta\vartheta
\let\epsilon\varepsilon

\newtheorem{theorem}{Theorem}
\newtheorem{lemma}[theorem]{Lemma}
\newtheorem{proposition}[theorem]{Proposition}
\newtheorem{corollary}[theorem]{Corollary}
\newenvironment{proof}[1][]%
{\setbox0=\hbox{#1}\paragraph{Proof\ifdim\wd0>0pt{ of #1:}\else{:}\fi}}%
{\hfill$\Box$\vskip1\baselineskip}

\newcommand{\Hypergeom}{\null_2F_1}
\newcommand{\close}[1]{\overline{#1}}
\newcommand{\Prob}{\mathop{\rm\bf\null P}\nolimits}
\newcommand{\Exp}{\mathop{\rm\bf\null E}\nolimits}
\newcommand{\C}{\mathbbm{C}}
\renewcommand{\H}{\mathbbm{H}}
\newcommand{\N}{\mathbbm{N}}
\newcommand{\R}{\mathbbm{R}}
\newcommand{\Z}{\mathbbm{Z}}
\renewcommand{\Re}{\mathop{\rm\null Re}\nolimits}
\renewcommand{\Im}{\mathop{\rm\null Im}\nolimits}
\newcommand{\ahat}{\hat{a}}
\newcommand\bhat{\hat{b}}
\newcommand\xhat{\hat{x}}
\newcommand{\im}{{\rm i}}
\newcommand{\floor}[1]{\left\lfloor{#1}\right\rfloor}
\newcommand{\ceil}[1]{\left\lceil{#1}\right\rceil}
\newcommand{\dif}[1][]{{\rm d}^{#1}}

\title{Reflected Brownian motion in generic triangles\\and wedges}
\author{Wouter Kager\thanks{Instituut voor Theoretische Fysica,
 Universiteit van Amsterdam, Valckenierstraat 65, 1018 XE Amsterdam,
 The Netherlands. E-mail: kager@science.uva.nl.}}

\begin{document}

\maketitle
\begin{abstract}
 Consider a generic triangle in the upper half of the complex plane with one
 side on the real line. This paper presents a tailored construction of a
 discrete random walk whose continuum limit is a Brownian motion in the
 triangle, reflected instantaneously on the left and right sides with
 constant reflection angles. Starting from the top of the triangle, it is
 evident from the construction that the reflected Brownian motion lands
 with the uniform distribution on the base. Combined with conformal
 invariance and the locality property, this uniform exit distribution allows
 us to compute distribution functions characterizing the hull generated by
 the reflected Brownian motion.
\end{abstract}

\section{Introduction and overview}
\label{sec:introduction}

\subsection{Motivation}
\label{ssec:motivation}

Reflected Brownian motions in a wedge with constant reflection angles on the
two sides were characterized by Varadhan and Williams in~\cite{varadhan:1985}.
Recent work by Lawler, Schramm and Werner~\cite{lsw:2003} on SLE establishes
a connection between one of these reflected Brownian motions (with reflection
angles of~$60^\circ$ with respect to the boundary), chordal SLE$_6$ and the
exploration process of critical percolation. They show that these three
processes generate the same hull, using an argument that we shall outline in
the following paragraph. The connection was carried even further by Julien
Dub\'edat~\cite{dubedat:2003,dubedat:2004}. He compared SLE$_6$ and the
aforementioned reflected Brownian motion in an equilateral triangle, started
from a given corner and conditioned to cross the triangle to a given point~$X$
on the opposite side. For these two processes, he showed that the conditional
probability that the side to the right of the starting point is the last side
visited before reaching~$X$, is the same. He also showed that, under the
assumption that SLE$_6$ is the scaling limit of the exploration process of
critical percolation, this conditional probability can be used to prove
Watts' formula for critical percolation (see~\cite{dubedat:2004}).

The argument used in~\cite{lsw:2003} to establish the connection between
the three models is as follows. Let~$(Z_t)$ be a stochastic process in an
equilateral triangle~$T$, started from a corner and stopped when it first
hits the opposite side~$S$. Define the \emph{hull}~$K_t$ as the compact set
of points in~$T$ disconnected from~$S$ by the trace of~$Z_t$ up to time~$t$.
Denote by~$X$ the hitting point of~$S$, and by~$\tau$ the hitting time. We
shall refer to the probability distribution of the point~$X$ as the ``exit
distribution''. For all three processes mentioned above, the exit distribution
is uniform on the side of the triangle. Together with conformal invariance
and the locality property (see Section~\ref{ssec:invariance} below
and~\cite{lsw:I}), shared by all three processes, this exit distribution
determines the law of the hull~$K_\tau$. The three processes therefore
generate the same hull.

This argument can be generalized to stochastic processes in arbitrary triangles
that are conformally invariant and have the locality property: if the exit
distributions of two such processes are the same, in particular if they are
both uniform, then the processes generate the same hull. The papers of Lawler,
Schramm, Werner~\cite{lsw:2003} and Dub\'edat~\cite{dubedat:2003} show that
in equilateral and isosceles triangles there exist reflected Brownian motions
that have uniform exit distributions. The purpose of this paper is to
generalize this result to arbitrary triangles, to review properties of these
reflected Brownian motions and to discuss distribution functions associated
with these processes and the hulls they generate.

\subsection{Notations and overview}
\label{ssec:notations}

To present an overview of the present paper, we first need to introduce
some notation. Given two angles $\alpha,\beta\in(0,\pi)$ such
that $\alpha+\beta<\pi$, we define the wedge $W_{\alpha,\beta}$ as the set
$\{z\in\C:\alpha-\pi < \arg z < -\beta\}$. We also define $T_{\alpha,\beta}$
as the triangle in the upper half of the complex plane such that one side
coincides with the interval~$(0,1)$, and the interior angles at the corners
$0$ and~$1$ are equal to $\alpha$ and~$\beta$, respectively. The third corner
is at $w_{\alpha,\beta}:=(\cos\alpha\sin\beta + \im\sin\alpha\sin\beta)/
\sin(\alpha+\beta)$. When $\alpha+\beta\geq\pi$, the domain $T_{\alpha,\beta}$
is similarly defined as the (unbounded) polygon having one side equal to the
interval~$(0,1)$, and interior angles $\alpha$ and~$\beta$ at the corners
$0$ and~$1$. We then identify the point at~$\infty$ with the ``third
corner''~$w_{\alpha,\beta}$.

Suppose now that $\theta_L$ and~$\theta_R$ are two angles of reflection on the
left and right sides of the wedge $W_{\alpha,\beta}$, respectively, measured
from the boundary with small angles denoting reflection away from the origin
($\theta_L,\theta_R\in(0,\pi)$). We shall use the abbreviation RBM$_{\theta_L,
\theta_R}$ to denote the corresponding reflected Brownian motion in the wedge
$W_{\alpha,\beta}$. For a characterization and properties of these RBMs, see
Varadhan and Williams~\cite{varadhan:1985} (note that here we use a different
convention for the reflection angles, namely, the angles $\theta_1$
and~$\theta_2$ of Varadhan and Williams correspond in our notation to the
angles $\theta_L-\pi/2$ and $\theta_R-\pi/2$).

The main goal of this paper is to show that in every wedge~$W_{\alpha,\beta}$
there is a unique RBM$_{\theta_L,\theta_R}$ with the following property:
started from the origin, the first hitting point of the RBM of any horizontal
line segment intersecting the wedge is uniformly distributed. This special
behaviour is obtained by taking the reflection angles equal to the angles of
the wedge, that is, $\theta_L=\alpha$ and $\theta_R=\beta$. Restricting
the wedge to a triangle we can reformulate this result as follows:

\begin{theorem}
 \label{the:main}
 Let $\alpha,\beta\in(0,\pi)$, $\alpha+\beta<\pi$, and let~$(Z_t:t\geq0)$
 be an RBM$_{\alpha,\beta}$ in the triangle $T_{\alpha,\beta}$ started
 from~$w_{\alpha,\beta}$ and stopped when it hits~$[0,1]$. Set $\tau:=
 \inf\{t>0:Z_t\in[0,1]\}$ and $X:=Z_\tau$. Then~$X$ is uniform in~$[0,1]$.
\end{theorem}

To prove this theorem we will cover the wedge~$W_{\alpha,\beta}$ with a
well-chosen lattice, and then define a random walk on this lattice. By
construction, this random walk will have the desired property that it
arrives on each horizontal row of vertices on the lattice with the uniform
distribution. Taking the scaling limit then yields the desired result. The
proof is split in two sections. In Section~\ref{sec:restricted} we consider
the easier case where the angles $\alpha$ and~$\beta$ are in the range
$(0,\pi/2]$ such that $\pi/2\leq\alpha+\beta<\pi$. Section~\ref{sec:generic}
treats the extension to arbitrary triangles, which is considerably more
involved.

Section~\ref{sec:properties} collects some properties of the two-parameter
family of RBMs. In particular, as was noted by Dub\'edat~\cite{dubedat:2003},
we can show from the discrete approximations that the imaginary parts of
the RBMs are essentially 3-dimensional Bessel processes. This allows us to
describe the time-reversals of the RBMs, and sheds some light on how the
RBM$_{\alpha,\beta}$ behaves in the domain~$T_{\alpha,\beta}$ for angles
$\alpha,\beta\in(0,\pi)$ such that $\alpha+\beta\geq\pi$. Finally, in
Section~\ref{sec:distributions} we compute several distribution functions
associated with the RBMs and the hulls they generate. This also reveals
intriguing connections between RBMs started from different corners of the
same triangle.

\section{RBM in a restricted geometry}
\label{sec:restricted}

Throughout this section, we assume that the angles $\alpha$ and~$\beta$ are
fixed and restricted to the range~$(0,\pi/2]$ such that $\pi/2\leq\alpha+
\beta<\pi$. We shall construct a random walk in the wedge~$W=W_{\alpha,
\beta}$ whose scaling limit is an RBM$_{\alpha,\beta}$ and whose horizontal
coordinate is uniformly distributed all the time, and thereby prove
Theorem~\ref{the:main} for these restricted values of $\alpha$ and~$\beta$.
The generalization to a generic triangle will be treated in
Section~\ref{sec:generic}.

\subsection{The lattice and further notations}
\label{ssec:lattice}

\begin{figure}
 \centering\includegraphics{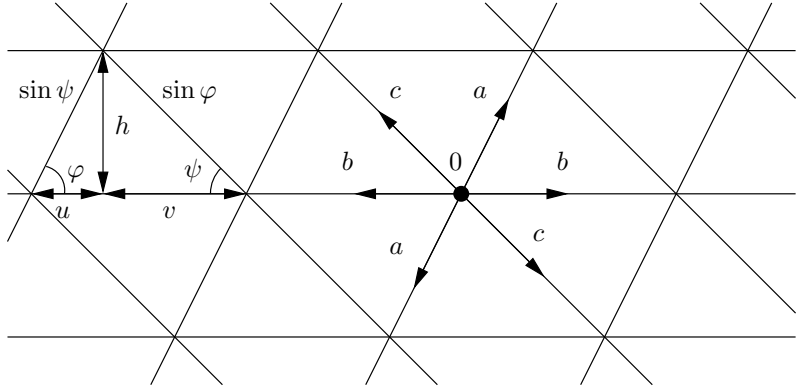}
 \caption{Picture of the lattice, showing the dimensions on the left, and
  the transition probabilities for a step of the random walk (from the
  origin in this picture) on the right.}
 \label{fig:Lattice}
\end{figure}

Throughout this paper we shall make use of a distorted triangular lattice,
defined as follows. Let $\phi$ and~$\psi$ be two angles in the range
$(0,\pi/2]$ such that $\pi/2\leq\phi+\psi<\pi$ (as we shall see later on, the
range for $\phi$ and~$\psi$ is chosen such that the transition probabilities
for our random walk are positive). We define $\Gamma_{\phi,\psi}$ as the
set of vertices $\{j\sin(\phi+\psi) - k\exp(\im\phi)\sin\psi : j,k\in\Z\}$
building the triangular lattice depicted in Figure~\ref{fig:Lattice}.
Throughout the paper we shall make use of the variables $u:=\cos\phi\sin\psi$,
$v:=\sin\phi\cos\psi$ and~$h:=\sin\phi\sin\psi$ to denote the lattice
dimensions. When we use these variables, the values of $\phi$ and~$\psi$
will always be clear from the context.

To introduce some further notation, let us first consider how one defines
a random walk $(X_n:n\geq0)$ on~$\Gamma_{\phi,\psi}$ that converges to
standard complex Brownian motion in the full plane, before we consider the
random walk in the wedge~$W$ in the following subsection. We set $X_0:=0$
and for each $n>0$, $X_n$ is chosen among the nearest-neighbours of~$X_{n-1}$
according to the probabilities $a$, $b$ and~$c$ as depicted in
Figure~\ref{fig:Lattice}. We may then write the position~$X_n$ of the random
walk as a sum of steps $S_1+S_2+\cdots+S_n$ where each step~$S_n= U_n +
\im V_n$ is a complex-valued random variable taking on the possible values
\begin{equation}
 S_n = \left\{\begin{array}{ll}
  \pm(u+\im h) & \mbox{with probability}~a; \\
  \pm(u+v) & \mbox{with probability}~b; \\
  \pm(v-\im h) & \mbox{with probability}~c.
  \end{array}\right.
\end{equation}

To obtain a two-dimensional Brownian motion as the scaling limit of the
random walk~$(X_n)$, it is sufficient that the covariance matrix of the
real and imaginary parts $U_n$ and~$V_n$ of each step is a multiple of
the identity. This gives two equations for the probabilities $a$, $b$
and~$c$:
\begin{equation}
 (a+b)\cot^2\phi + 2b\cot\phi\cot\psi + (b+c)\cot^2\psi = a+c, \\
\end{equation}
\begin{equation}
 a\cot\phi - c\cot\psi = 0,
\end{equation}
where $\cot x=1/\tan x$. The probabilities $a$, $b$ and~$c$ can be determined
from these equations, yielding
\begin{eqnarray}
 a &=& \lambda\cot\psi(\cot\phi+\cot\psi), \label{equ:a} \\
 b &=& \lambda(1-\cot\phi\cot\psi), \label{equ:b} \\
 c &=& \lambda\cot\phi(\cot\phi+\cot\psi),  \label{equ:c}
\end{eqnarray}
where $\lambda=\frac{1}{2}[\cot\phi(\cot\phi+\cot\psi)+\sin^{-2}\psi]^{-1}$
is the normalization constant. One may verify that $\phi$ and~$\psi$ must
satisfy $\pi/2\leq\phi+\psi<\pi$ to make all three probabilities nonnegative.

We conclude this subsection with a short discussion of how one obtains the
scaling limit of the random walk~$(X_n)$. To do so, for every natural
number~$N>0$ one may define the continuous-time, complex-valued stochastic
process~$(Z^{(N)}_t:t\geq0)$ as the linear interpolation of the process
\begin{equation}
 Y^{(N)}_t = \frac{1}{N\sigma}X_{\floor{N^2 t}}
\end{equation}
making jumps at the times $\{k/N^2:k=1,2,\ldots\}$. Here, $\sigma^2$ is
the variance of the real and imaginary parts of the steps~$S_n$, that is,
$\sigma^2=\Exp[U_n^2]=\Exp[V_n^2]$. It is then standard that $Z^{(N)}_t$
converges weakly to a complex Brownian motion in the full complex plane~$\C$
when $N\to\infty$ (topological aspects are as described in the Introduction
of Varadhan and Williams~\cite{varadhan:1985}).

\subsection{Reflected random walk in a wedge}
\label{ssec:wedge}

\begin{figure}
 \centering\includegraphics{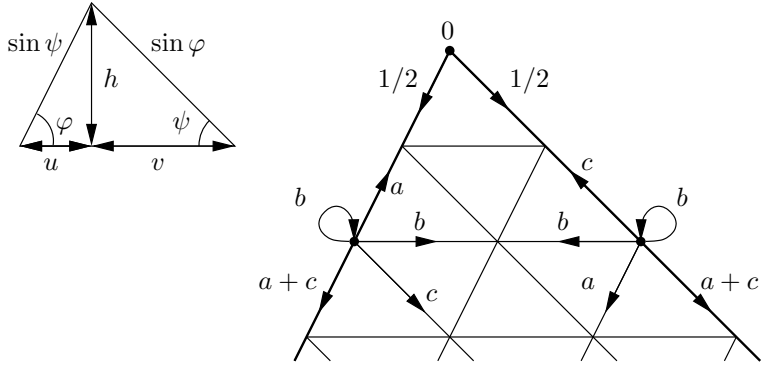}
 \caption{Definition of the reflected random walk in a wedge.}
 \label{fig:WedgeSimple}
\end{figure}

We now return to the wedge $W=W_{\alpha,\beta}$ with angles $\alpha$
and~$\beta$ fixed in the range~$(0,\pi/2]$ such that $\pi/2\leq\alpha+\beta
<\pi$. Clearly, this wedge is covered nicely with vertices of the
lattice~$\Gamma_{\phi,\psi}$ when we set $\phi:=\alpha$ and~$\psi:=\beta$,
see Figure~\ref{fig:WedgeSimple}. For the duration of this section we
consider these values of $\phi$ and~$\psi$ to be fixed. Later, when we
generalize to arbitrary triangles, the relation between $\phi,\psi$ and
$\alpha,\beta$ will not be so simple, which is why we already reserve the
symbols $\phi,\psi$ to denote the angles of the lattice.

We shall denote by~$G=G_{\alpha,\beta}$ the set of vertices obtained by
taking the intersection of~$\Gamma=\Gamma_{\phi,\psi}$ with~$\close{W}$.
We shall call the vertices of~$G$ having six nearest neighbours along the
lattice directions \emph{interior vertices}. The origin will be called the
\emph{apex} of~$G$, and the remaining vertices will be referred to as the
\emph{boundary vertices}. The set of boundary vertices may be further
subdivided into \emph{left boundary vertices} and \emph{right boundary
vertices}, with the obvious interpretation.

Given a vertex~$x$ of~$G$, a reflected random walk $(X^x_n:n\geq0)$ on~$G$ is
defined as follows. We set $X^x_0:=x$ and for each~$n>0$, if $X^x_{n-1}$ is
an interior vertex, then $X^x_n$ is chosen among the six nearest-neighbours
of~$X^x_{n-1}$ according to the probabilities $a$, $b$ and~$c$ as before.
This guarantees that the scaling limit of the random walk is Brownian
motion in the interior of the wedge. It remains to specify the transition
probabilities for the random walk from the boundary vertices and the apex.

In this paper, we restrict ourselves to the case where the transition
probabilities on all left boundary vertices are the same, and likewise for
the right boundary vertices. We write $\Exp_z[S_1]$ for the expected value
of the first step of the random walk started from the vertex~$z$, and
assume that it is nonzero at all the boundary vertices and the apex. Then,
as we shall prove in Section~\ref{ssec:scalinglimit}, the random walk
converges to a reflected Brownian motion. Moreover, the directions of
reflection on the sides of the wedge are given by the directions of
$\Exp_z[S_1]$ at the left and right boundary vertices. It follows that by
playing with the transition probabilities from the boundary vertices we
can obtain different RBMs in the scaling limit.

There is, however, only one choice of transition probabilities for which
the random walk has the following special property: started from the origin,
the random walk first arrives on any row of the lattice with the uniform
distribution on the vertices of that row. We will refer to this special case
as the \emph{uniform (random) walk}. To derive its transition probabilities
one proceeds as follows. In a picture where we represent the steps of the
walk by arrows, we have to make sure that every vertex in a given row has
two incoming arrows with probability~$b$ from vertices in the same row, and
two incoming arrows with probabilities $a$ and~$c$ from vertices in the rows
above and below. This completely determines the transition probabilities
from the boundary vertices, see Figure~\ref{fig:WedgeSimple} for a picture
of the solution.

In formula, if~$x$ is a left boundary vertex, then we have the following
transition probabilities for the uniform walk:
\begin{eqnarray}
 && p[x,x] = p[x,x +(u+v)] = b, \\
 && p[x,x + u + \im h] = a, \\
 && p[x,x - u - \im h] = a+c, \\
 && p[x,x + v - \im h] = c,
\end{eqnarray}
and the transition probabilities from the right boundary vertices are defined
symmetrically, as shown in Figure~\ref{fig:WedgeSimple}. At the apex we
simply choose the transition probabilities to each the two vertices directly
below the apex equal to~$1/2$.

It will now be convenient to decompose the position~$X^x_n$ at each step of
the uniform walk as $J_n(u+v) - K_n(u+\im h)$. Then $K_n$ is a nonnegative
integer denoting a row of vertices on the lattice, and $J_n$ is a nonnegative
integer denoting the position on the $K_n$th row. Observe that there are a
total of $N(k)=k+1$ vertices on the $k$th row, so that $J_n$ ranges from $0$
to $N(K_n)-1$. Henceforth, we shall always adopt this convention of numbering
rows on the lattice in top-down order, and vertices on each row from left to
right.

By construction, the uniform random walk~$(X^0_n)$ started from the origin
has the property that at every time $n\geq0$, the position of the walker is
uniformly distributed on the rows of the lattice. More precisely, for the
uniform walk~$(X^0_n)$ the following lemma holds:

\begin{lemma}
 \label{lem:uniform}
 For all~$n\in\N$, if $k_0,k_1,\ldots,k_n$ is a sequence of natural numbers
 such that $k_0=0$ and $|k_m-k_{m-1}|\leq1$ for all $m=1,2,\ldots,n$, then
 for each $j=0,1,\ldots,N(k_n)-1$,
 \[ \Prob[J_n=j\mid K_0=k_0,\ldots,K_n=k_n]
     = \Prob[J_n=j\mid K_n=k_n] = \frac{1}{N(k_n)}.
 \]
 In particular, if $k\in\N$ and $T:=\min\{n\geq0:K_n\in k\}$ is the first
 time at which the walk visits row~$k$, then for each $j=0,1,\ldots,N(k)-1$,
 \[ \Prob[J_T=j] = \frac{1}{N(k)}.
 \]
\end{lemma}

\begin{proof}
 The first claim of the lemma is proved by induction. For $n=0,1$ the claim
 is trivial. For~$n>1$,
 \begin{eqnarray}
  \lefteqn{\Prob[J_n = j \mid K_n=k_n]\phantom{\int_I}} \nonumber\\
  && = \frac{\sum\limits_{k_{n-1}} \Prob[J_n=j, K_n=k_n\mid K_{n-1}=k_{n-1}]\,
           \Prob[K_{n-1}=k_{n-1}]}
          {\sum\limits_{k_{n-1}} \Prob[K_n=k_n\mid K_{n-1}=k_{n-1}]\,
           \Prob[K_{n-1}=k_{n-1}]}.
  \label{equ:uniform}
 \end{eqnarray}
 Using the induction hypothesis, it is easy to compute the conditional
 probabilities in the numerator and denominator of this expression from the
 transition probabilities given earlier. We first assume~$k_n>1$. Then for
 $k_{n-1}=k_n$ these conditional probabilities are $2b/N(k_n)$ and $2b$,
 respectively. For $k_{n-1}=k_n\pm1$ they are $(a+c)/N(k_{n-1})$ and
 $(a+c)N(k_n)/N(k_{n-1})$, and for all other~$k_{n-1}$ they vanish. When
 $k_n\leq1$ the contributions from the apex take a different form, but the
 computation is equally straightforward. Now observe that the quotient in
 Equation~(\ref{equ:uniform}) has the same value if the event $\{K_{n-1}=
 k_{n-1}\}$ is replaced by $\{K_0=k_0,\ldots,K_{n-1}=k_{n-1}\}$. The first
 claim of the lemma follows. The second claim of the lemma is then a
 straightforward consequence.
\end{proof}

To prove convergence of the uniform walk to a reflected Brownian motion, we
shall need an estimate on the local time spent by the walk on the boundary
vertices of~$G$. This estimate is provided by Lemma~\ref{lem:localtime}
below. The proof of Lemma~\ref{lem:localtime} is postponed to
Section~\ref{ssec:intertwiningdiscrete}, since the proof will give us some
intermediate results that we will need only in Section~\ref{sec:properties}.

\begin{lemma}
 Consider the walk~$(X^0_n)$, killed when it reaches row~$M>0$, and let~$x$
 be a fixed boundary vertex of~$G$. Then the expected number of visits
 of~$(X^0_n)$ to~$x$ before it is killed is smaller than $(a+c)^{-1}$.
 \label{lem:localtime}
\end{lemma}

\subsection{Identification of the scaling limit}
\label{ssec:scalinglimit}

Given a point~$x$ in~$W$, the scaling limit of the reflected random walk
started from a vertex near~$x$ is obtained in a similar way as described in
Section~\ref{ssec:lattice}. That is, if $(x_N)$ is a sequence in~$G$ such
that $|x_N-N\sigma x|\leq1/2$ for all $N\in\N$, then for every natural number
$N>0$ we define the continuous complex-valued stochastic
process~$(Z^{(N)}_t : t\geq0)$ as the linear interpolation of the process
\begin{equation}
 Y^{(N)}_t = \frac{1}{N\sigma}X^{x_N}_{\floor{N^2t}}
\end{equation}
making jumps at the times $\{k/N^2:k=1,2,\ldots\}$. Here, $\sigma^2$ is again
the variance of the real and imaginary parts of the steps $S_n=X_n-X_{n-1}$
of the random walk in the interior of the wedge. Then it is clear that the
random walk converges to a complex Brownian motion in the interior, but the
behaviour on the boundary is non-trivial.

It is explained in Dub\'edat~\cite{dubedat:2003} how one proves that the
scaling limit of the random walk~$(X^x_n)$ is in fact reflected Brownian
motion with constant reflection angles on the two sides, started from~$x$.
The core of the argument identifies reflected Brownian motion as the only
possible weak limit~$(Z^x_t)$ of~$(Z^{(N)}_t)$. We repeat this part of the
argument below, mainly to allow us in Section~\ref{sec:generic} to point
out some technical issues that need to be taken care of when we generalize
to arbitrary triangles.

We use the submartingale characterization of reflected Brownian motion
in~$W$ (see~\cite{varadhan:1985}, Theorem~2.1), which states the following.
Let~$C^2_b(W)$ be the set of bounded continuous real-valued functions on~$W$
that are twice continuously differentiable with bounded derivatives. Assume
that the two reflection angles $\theta_L$ and~$\theta_R$ on the sides of~$W$
are given. Then the RBM$_{\theta_L,\theta_R}$ in~$W$ started from $x\in W$ is
the unique continuous strong Markov process~$(Z^x_t)$ in~$W$ started from~$x$
such that for any $f\in C^2_b(W)$ with nonnegative derivatives on the
boundary in the directions of reflection, the process
\begin{equation}
 f(Z^x_t)-\frac{1}{2}\int_0^t \Delta f(Z^x_s)\,\dif{s},
\end{equation}
where~$\Delta$ is the Laplace operator, is a submartingale.

To prove convergence of our random walk to an RBM, it is therefore sufficient
to show that there are two angles $\theta_L$ and~$\theta_R$ such that if~$f$
is a function as described above, then for all $0\leq s<t$
\begin{equation}
 \liminf_{N\to\infty} \Exp^{(N)}_x\left[ f(Z_t)-f(Z_s)
  - \frac{1}{2}\int_s^t \Delta f(Z_u)\,\dif{u} \right] \geq 0.
 \label{equ:submartingale}
\end{equation}
Here, $\Exp^{(N)}_x$ denotes the expectation operator for the $N$th
approximate process started near~$x$, and we have dropped the $(N)$
superscript on $(Z_t)$ to simplify the notation. Note that the conditioning
in the submartingale property is taken care of because the starting point~$x$ 
is arbitrary. In fact, it is sufficient to verify~(\ref{equ:submartingale})
up to the stopping time $\tau:=\inf\{t\geq0:\Im Z_t\leq-M\}$ for some large
number~$M$. We will make use of this to have uniform bounds on the error
terms in our discrete approximation.

We now show~(\ref{equ:submartingale}). Let~$f$ be a function in~$C^2_b(W)$,
and write $u_k=k/N^2$ for the jump times of~$(Z_t)$. Then by Taylor's theorem,
there exist (random) times~$T_k$ between $u_k$ and~$u_{k+1}$ such that
\begin{eqnarray}
 \lefteqn{ \hskip-1em f(Z_{u_{k+1}}) - f(Z_{u_k}) =
  \frac{1}{N\sigma}\Big( \partial_x f(Z_{u_k})U_{k+1}
  + \partial_y f(Z_{u_k})V_{k+1} \Big) + \null} \nonumber\\
  && \hskip-2em
   \frac{1}{2N^2\sigma^2}\Big( \partial^2_x f(Z_{T_k})U^2_{k+1}+
   \partial^2_y f(Z_{T_k})V^2_{k+1}
   + 2\partial_x\partial_y f(Z_{T_k})U_{k+1}V_{k+1} \Big),
\end{eqnarray}
where $U_k$ and~$V_k$ are the real and imaginary parts of the $k$th step
of the underlying random walk, as before.

When we now apply the expectation operator, we distinguish between the
behaviour on the boundary and in the interior (we may ignore what happens
at the apex because of the negligible time spent there, see below). If~$B$
denotes the set of boundary vertices of~$G$, then we may write
\begin{eqnarray}
 && \hskip-2.5em\Exp^{(N)}_x\left[f(Z_{u_{k+1}}) - f(Z_{u_k})\right] =
  \Exp^{(N)}_x\Bigg[ \frac{1}{2N^2}\Delta f(Z_{T_k}) + \null
  \nonumber\\
 && \hskip-2em
    \left( \frac{1}{N\sigma}\partial_x f(Z_{u_k})U_{k+1}
  + \frac{1}{N\sigma}\partial_y f(Z_{u_k})V_{k+1} + O(N^{-2}) \right)
    \mathbbm{1}_{\{N\sigma Z_{u_k}\in B\}} \Bigg],
\end{eqnarray}
where the $O(N^{-2})$ error term is uniform in $\{z\in W:\Im z > -M\}$. Summing
over~$k$ from $\floor{N^2(s\land\tau)}$ to $\floor{N^2(t\land\tau)}$ yields
\begin{eqnarray}
 && \hskip-3.5em
 \Exp^{(N)}_x\left[ f(Z_{t\land\tau})-f(Z_{s\land\tau}) \right] =
 \Exp^{(N)}_x\left[
  \int_{s\land\tau}^{t\land\tau}\frac{1}{2}\Delta f(Z_u)\dif{u} \right]
  + o(1) + \null\nonumber\\
 && \hskip-3.5em
 \sum_{z\in B}\Exp^{(N)}_x[L_z]\left\{
  \frac{1}{N\sigma}\partial_x f(\frac{z}{N\sigma})\Exp_z[U_1]
  + \frac{1}{N\sigma}\partial_y f(\frac{z}{N\sigma})\Exp_z[V_1]
  + O(N^{-2})\right\}.
 \label{equ:braces}
\end{eqnarray}
In the last line, $L_z$ is the local time at~$z$ (the number of jumps of the
random walk to~$z$) between times $s\land\tau$ and $t\land\tau$, and $\Exp_z$
is expectation with respect to the underlying random walk started from~$z$.
Lemma~\ref{lem:localtime} shows that when the starting point of the walk
is~$0$, the expected local time spent on the boundary vertices up to
time~$\tau$ is of order~$O(N)$. Then one can use the Markov property to see
that this is enough for us to ignore the $O(N^{-2})$ term in the limit.

Observe that the remaining term in braces is just the derivative
of~$f$ at the boundary point~$z/N\sigma$ along the direction of $\Exp_z[S_1]$.
Thus, if the function~$f$ has nonnegative derivatives along this direction on
the boundary of~$W$, then the term in braces in equation~(\ref{equ:braces})
is nonnegative, and the desired result~(\ref{equ:submartingale}) is obtained.
This proves that the scaling limit of the reflected random walk is a
reflected Brownian motion, and that the angles of reflection with respect to
the left and right sides are given by
\begin{equation}
 \theta_L+\alpha=\arg\left(-\Exp_z[S_1]\right)
 \qquad\mbox{and}\qquad
 \theta_R+\beta=-\arg\left(\Exp_y[S_1]\right),
\end{equation}
where~$z$ is any left boundary vertex and~$y$ any right boundary vertex. It
is clear from the computation that this result is valid generally for any
choice of transition probabilities on the left and right sides. We shall
however focus again on the special case of the uniform walk introduced in
Section~\ref{ssec:wedge}.

\begin{figure}
 \centering\includegraphics{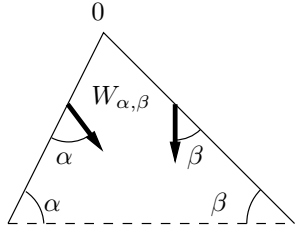}
 \caption{The thick arrows in this figure represent the directions of
  reflection of the reflected Brownian motion described in the text. The
  reflection angles are such that the Brownian motion will hit the dashed
  line with uniform distribution.}
 \label{fig:Angles}
\end{figure}

The computation of the reflection angles for the uniform walk is
straightforward. By symmetry it is enough to compute only~$\theta_R$. From
the transition probabilities one may verify that
\begin{equation}
 \cot(\theta_R+\beta)
  = -\frac{\Exp_z[\Re S_1]}{\Exp_z[\Im S_1]}
  = \frac{(a-b)\cot\psi-(a+b)\cot\phi}{2a}.
\end{equation}
Substituting equations~(\ref{equ:a})--(\ref{equ:c}) then yields
\begin{equation}
 \cot(\theta_R+\beta)
  = \frac{\cot^2\psi-1}{2\cot\psi}=\cot(2\psi)=\cot(2\beta).
\end{equation}
Thus, the angle of reflection of the Brownian motion with respect to the
right side is simply $\theta_R=\beta$. Similarly, the angle of reflection
on the left side is given by $\theta_L=\alpha$.

We conclude that the scaling limit of the uniform random walk defined in
Section~\ref{ssec:wedge} is a reflected Brownian motion with fixed reflection
angles on the two sides of the wedge. The angles of reflection are $\alpha$
and~$\beta$ with respect to the left and right sides, respectively, as
illustrated in Figure~\ref{fig:Angles}. It follows from Lemma~\ref{lem:uniform}
that the RBM$_{\alpha,\beta}$ has the special property that in the
wedge~$W_{\alpha,\beta}$, the RBM first arrives on any horizontal
cross-section of the wedge with the uniform distribution. By
a simple translation it follows that in the triangle~$T_{\alpha,\beta}$,
the RBM$_{\alpha,\beta}$ started from the top~$w_{\alpha,\beta}$ will land
on~$[0,1]$ with the uniform distribution. This completes the proof of
Theorem~\ref{the:main} for angles $\alpha,\beta\in(0,\pi/2]$ that satisfy
$\pi/2\leq\alpha+\beta<\pi$. The case of arbitrary angles will be considered
in the following section.

\section{RBM in a generic triangle}
\label{sec:generic}

In the previous section we proved Theorem~\ref{the:main} for angles $\alpha,
\beta\in(0,\pi/2]$ such that $\pi/2\leq\alpha+\beta<\pi$. Here we generalize
to generic triangles. That is, we now choose the angles $\alpha$ and~$\beta$
arbitrarily in the range~$(0,\pi)$ such that $\alpha+\beta<\pi$. These values
of $\alpha$ and~$\beta$ are assumed fixed for the remainder of this section.
As in Section~\ref{sec:restricted}, we will define a uniform walk in the
wedge~$W_{\alpha,\beta}$ and identify the scaling limit as a reflected
Brownian motion.

\subsection{Choice of the lattice}
\label{ssec:lattice_gen}

The first thing that we have to do is to find a lattice covering the wedge
$W=W_{\alpha,\beta}$ in a nice way. We will make use of the triangular
lattices $\Gamma_{\phi,\psi}$ introduced in the previous section. Consider
such a lattice~$\Gamma_{\phi,\psi}$, and choose two vertices $x_L$ and~$x_R$
on the first row of the lattice such that~$x_R$ is to the right of~$x_L$.
Then the two half-lines $\{tx_L:t\geq 0\}$ and $\{tx_R:t\geq0\}$ define the
left and right sides of a wedge that is covered nicely by vertices
of~$\Gamma_{\phi,\psi}$. We will show below that for any given wedge~$W$
there is a choice of the lattice angles $\phi,\psi$ and the two points
$x_L$ and~$x_R$ such that the wedge thus defined coincides with~$W$.

\begin{figure}
 \centering\includegraphics{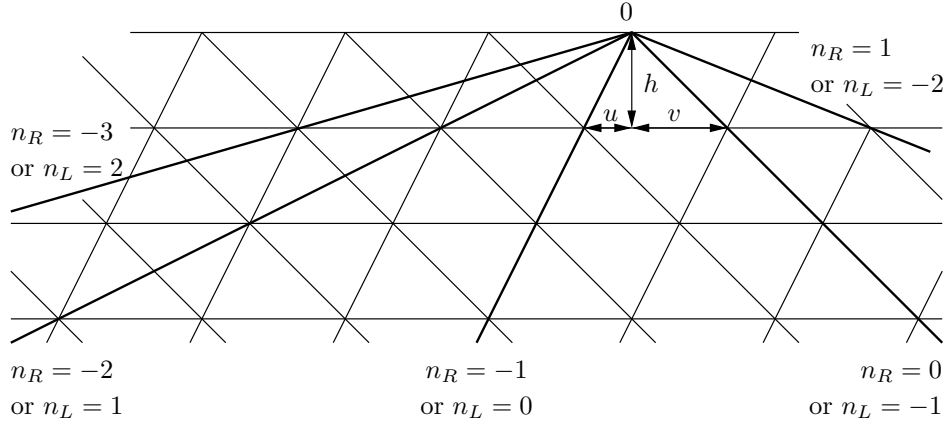}
 \caption{Different wedges can be covered by the same triangular lattice,
  by changing the directions of the two sides as shown. These directions
  can be expressed in terms of integers $n_L$ and~$n_R$, as explained in
  the text. Each thick line in the figure can represent either the right
  side of a wedge (for which the corresponding value of~$n_R$ is given in
  the figure), or the left side of a wedge (for which the corresponding
  value of~$n_L$ is given).}
 \label{fig:ExtraVertices}
\end{figure}

We start by introducing some notation. Given the lattice~$\Gamma_{\phi,\psi}$
and the two vertices $x_L$ and~$x_R$, we may introduce two integers $n_L$
and~$n_R$ that count the positions of $x_L$ and~$x_R$ on the first row of
the lattice. More precisely, the integers $n_L$ and~$n_R$ are defined such
that $x_L=-u-n_L(u+v)-\im h$ and $x_R=v+n_R(u+v)-\im h$ (for an illustration,
see Figure~\ref{fig:ExtraVertices}). Conversely, given two integers $n_L$
and~$n_R$ and the lattice~$\Gamma_{\phi,\psi}$, the vertices $x_L$ and~$x_R$
are fixed by these equations. Observe that $n_L$ and~$n_R$ must satisfy
$n_L+n_R\geq0$ to make sure that $x_R$ lies to the right of~$x_L$.
Figure~\ref{fig:ExtraVertices} shows the wedges one obtains for different
choices of $n_L$ and~$n_R$ on a given lattice.

The main claim of this subsection is that for any  choice of the angles
$\alpha$ and~$\beta$, there is a choice of integers $n_L,n_R$ and of the
lattice angles $\phi,\psi$ such that the wedge one obtains as described
above coincides with the wedge~$W_{\alpha,\beta}$. This results is a direct
consequence of Lemma~\ref{lem:geometry} below. The proof of the lemma,
which will give explicit formulas for $n_L,n_R$ and $\phi, \psi$ in terms
of the angles $\alpha$ and~$\beta$, is postponed to the end of this
subsection.

\begin{lemma}
 \label{lem:geometry}
 Let $\alpha,\beta\in(0,\pi)$ such that $\alpha+\beta<\pi$. Then
 there is a choice of (possibly negative) integers $n_L$ and~$n_R$
 with $n_L + n_R \geq 0$, and angles $\phi,\psi\in(0,\pi/2]$ with
 $\pi/2\leq\phi+\psi<\pi$, such that
 \begin{eqnarray}
  \cot\alpha &=& n_L(\cot\phi+\cot\psi) + \cot\phi;
   \label{equ:alpha} \\
  \cot\beta  &=& n_R (\cot\phi+\cot\psi) + \cot\psi.
   \label{equ:beta}
 \end{eqnarray}
\end{lemma}

We remind the reader that in Section~\ref{sec:restricted} we considered a
random walk on a wedge with angles $\alpha,\beta\in(0,\pi/2]$ such that
$\pi/2<\alpha+\beta<\pi$. The situation described there corresponds to the
special choice of $n_L=n_R=0$, $\phi=\alpha$ and $\psi=\beta$ in
Lemma~\ref{lem:geometry}. Thus, the random walk construction of
Section~\ref{sec:restricted} is a special case of the more general
construction we are considering in this section. We would also like to
remark at this point that the choice of $n_L,n_R$ and $\phi,\psi$ in
Lemma~\ref{lem:geometry} is not unique in general. For instance, we will
see in the proof of Lemma~\ref{lem:geometry} that it is always possible
to choose the angles $\phi,\psi$ in the range $[\pi/4,\pi/2]$. Thus, there
are angles $\alpha,\beta$ such that the choices made in the proof of the
lemma and in the construction of Section~\ref{sec:restricted} are different.

From now on, we will assume that the values of $\phi$ and~$\psi$ are fixed as
in Lemma~\ref{lem:geometry}, so that the lattice $\Gamma=\Gamma_{\phi,\psi}$
is fixed. As we explained above, the lattice provides a nice covering of the
wedge~$W=W_{\alpha,\beta}$. Following the conventions introduced in
Section~\ref{sec:restricted}, we shall denote by~$G=G_{\alpha,\beta}$ the set
of vertices obtained by taking the intersection of~$\Gamma$ with~$\close{W}$.
We again subdivide the set~$G$ into the \emph{apex}, \emph{interior vertices}
and \emph{left} and \emph{right boundary vertices}. We now conclude this
subsection with the proof of Lemma~\ref{lem:geometry}.

\begin{proof}[Lemma~\ref{lem:geometry}]
Assume first that both $\alpha$ and~$\beta$ are smaller than~$\frac{\pi}{2}$.
Then we can take
\begin{equation}
 n_L = \ceil{\cot\alpha} - 1,\quad
 n_R = \ceil{\cot\beta} - 1,
\end{equation}
and solve equations (\ref{equ:alpha}) and~(\ref{equ:beta}) for $\phi$
and~$\psi$ to obtain
\begin{eqnarray}
 \cot\phi &=& \frac{n_R+1}{n_L+n_R+1}\cot\alpha
  - \frac{n_L}{n_L+n_R+1}\cot\beta; \label{equ:phi} \\
 \cot\psi &=& \frac{n_L+1}{n_L+n_R+1}\cot\beta
  - \frac{n_R}{n_L+n_R+1}\cot\alpha. \label{equ:psi}
\end{eqnarray}
Observe that since $n_L<\cot\alpha\leq n_L+1$ and
$n_R<\cot\beta\leq n_R+1$, the angles $\phi$ and~$\psi$ are in the
range~$[\pi/4,\pi/2]$.

It remains to consider the case when either $\alpha$ or~$\beta$ is at
least~$\pi/2$, and by symmetry it suffices to assume $\alpha\geq\pi/2$.
Then, we can for instance set $k=\ceil{\cot\alpha+\cot\beta}$, and let~$l$
be the smallest positive integer such that $l/(k+l)>-\cot\alpha/\cot\beta$.
We then set
\begin{equation}
 n_L := -l, \quad n_R := k+l-1,
\end{equation}
and the angles $\phi$ and~$\psi$ are given by the equations (\ref{equ:phi})
and~(\ref{equ:psi}) as before. From the inequalities
\begin{equation}
 \frac{l}{k+l} > -\frac{\cot\alpha}{\cot\beta} \geq \frac{l-1}{k+l-1},
\end{equation}
plus the fact that $k\geq\cot\alpha+\cot\beta$, it follows that $\phi$
and~$\psi$ are in the range~$[\pi/4,\pi/2]$. This completes the proof.
\end{proof}

\subsection{Discussion of the random walk construction}
\label{ssec:genericdiscussion}

Now that we have chosen the lattice on the wedge~$W$, the next task is
to define a random walk~$(X^x_n)$ on~$G$ whose scaling limit is reflected
Brownian motion. As in Section~\ref{ssec:wedge}, we will focus on the special
case of the uniform walk, i.e., the random walk on the lattice that stays
uniform on the rows all the time. Earlier we explained how one defines this
random walk in the case $n_L=n_R=0$. In this subsection we will describe
what one has to do to generalize to other wedges, i.e., to the case where
$n_L$ or~$n_R$ or both are nonzero. In the following subsection we will
then spell out the transition probabilities of the uniform walk for this
general case.

It is clear that we should define the transition probabilities from the
interior vertices in the same way as before. This will guarantee that the
random walk will converge to Brownian motion in the interior of the wedge.
The nontrivial task is to define the transition probabilities from the
boundary vertices and at the top of the wedge. Once these are defined, we
will use the strategy of Section~\ref{ssec:scalinglimit} to identify the
scaling limit as a reflected Brownian motion. Let us therefore reconsider
the arguments used there to identify the scaling limit, to see if they
still apply to general wedges.

Looking back at Section~\ref{ssec:scalinglimit} we see that it is still
sufficient to verify Equation~(\ref{equ:submartingale}) for appropriate
functions~$f$. Moreover, Equation~(\ref{equ:braces}) is still valid after
making a Taylor expansion. However, we emphasize that in the previous
section the sum over the boundary vertices picked up contributions from
only one left boundary vertex and one right boundary vertex on every row
of the lattice. This is no longer the case in the general situation, as
we shall now discuss.

Remember that the boundary vertices are defined as those vertices having less
than six nearest-neighbours in~$G$. Now consult Figure~\ref{fig:ExtraVertices}.
Then we see that for a given value of~$n_L$, the number of left boundary
vertices on each row of the lattice is fixed and equals $N_L=|n_L|+
1_{\{n_L\geq0\}}$. Likewise, for given~$n_R$ the number of right boundary
vertices on every row is $N_R=|n_R|+1_{\{n_R\geq0\}}$. (An observant reader
will notice that if $n_L$ or~$n_R$ is negative, then on the first few rows
of the lattice the number of boundary vertices is less than~$N_L+N_R$. This
is a complication that we will deal with later on.)

We conclude that the sum in Equation~(\ref{equ:braces}) picks up~$N_L$
contributions from the left boundary vertices in every row. Each of these
contributions is proportional to the derivative of~$f$ near the boundary in
the direction of $\Exp_z[S_1]$. We have to guarantee that these contributions
add up to something positive if~$f$ has positive derivative along some (yet
to be determined) fixed direction on the left side of the wedge. To do so,
it is clearly sufficient to impose that $\Exp_z[S_1]$ is the same at all
the left boundary vertices, and this is what we shall do. The corresponding
condition will be imposed at the right boundary vertices.

It is clear from the discussion in Section~\ref{ssec:scalinglimit} that
the scaling limit of the random walk will be a reflected Brownian motion.
Moreover, the directions of reflection on the two sides of~$W$ will be given
by the values of $\Exp_z[S_1]$ at the left and right boundary vertices. Let
us remark that the condition that $\Exp_z[S_1]$ is the same at all the boundary
vertices may be more than necessary, but by imposing this condition we avoid
having to estimate the differences in expected local time at different
boundary vertices.

In summary, we will look for transition probabilities from the boundary
vertices that satisfy the following conditions:
\begin{enumerate}
\item The summed transition probability to a given vertex from vertices in
the row above is $a+c$, and likewise from vertices in the row below.
\label{conditions}
\item The summed transition probability to a given vertex from vertices in
the same row is $2b$.
\item The expected value of the first step of the walk from every left
boundary vertex is the same, and likewise for the right boundary vertices.
\end{enumerate}
Observe that condition~1 introduces an up-down symmetry which is not inherent
in the geometry of the problem, but will be of importance later. Explicit
expressions for all the transition probabilities of the uniform walk will be
given in the following subsection. This will show that for any~$n_L,n_R$ it
is possible to choose the transition probabilities such that they satisfy
conditions~1--3.

To conclude this subsection, we point out one further complication that
we already commented on earlier. This complication arises when either $n_L$
or~$n_R$ is negative, because then the first few rows of the lattice will
contain less than $N_L+N_R$ vertices in total. Hence we will have to specify
separately the transition probabilities on these first rows of the lattice.
As we shall discuss below, this complication is quite easy to handle.

Consider once again Figure~\ref{fig:ExtraVertices} and recall the definition
of the integers $n_L$ and~$n_R$ in Section~\ref{ssec:lattice_gen}, illustrated
in the figure. Then one notes that each row of the lattice contains exactly
$n_L+n_R+1$ vertices more than the row above. In other words, the total
number of vertices on row~$k$ of the lattice is $N(k)=(n_L+n_R+1)k+1$. From
this one can compute that if either $n_L$ or~$n_R$ is negative, then the first
row of the lattice that contains at least all of the $N_L+N_R$ boundary
vertices is row $k_0$, where
\begin{equation}
 k_0:=\ceil{\frac{|n_L-n_R|}{n_L+n_R+1}}.
\end{equation}
For example, for $n_L=2$ and~$n_R=-1$ the numbers of left and right boundary
vertices are $N_L=3$ and~$N_R=1$, respectively. However, the first two rows
of~$G$ contain only $N(0)=1$ and~$N(1)=3$ vertices in total, so that the
first row of the lattice that contains at least $N_L+N_R$ vertices is row~$2$.

The remaining problem is then to define the transition probabilities from the
vertices in the first $k_0$~rows of the lattice. It is not very difficult to
choose the transition probabilities on these first few rows such that the
walk will trivially stay uniform on the rows of the lattice. In fact, we can
choose the transition probabilities such that conditions 1 and~2 above are
also satisfied at rows $1$ up to~$k_0-1$, as we shall see below. This
guarantees that the walk does indeed stay uniform on rows. We can not make
condition~3 hold on the first $k_0$~rows. This is no problem, since this
condition was only introduced to prove convergence to reflected Brownian
motion, as we discussed above, and the local time spent at the first~$k_0$
rows is negligible.

\subsection{Transition probabilities and scaling limit}
\label{ssec:walkgeneric}

In the previous subsection we described in words how one should define the
transition probabilities of the uniform walk in the generic wedge~$W=
W_{\alpha,\beta}$. Here we shall complete the uniform walk construction by
providing explicit expressions for all of the transition probabilities. We
then complete the proof of Theorem~\ref{the:main} by computing the directions
of reflection for the RBM obtained in the scaling limit.

To get us started, given a vertex~$x$ of~$G$ we set $X^x_0=x$. In
the interior of the wedge we define the transition probabilities as before.
That is, when~$X^x_n$ is an interior vertex, $X^x_{n+1}$ is to be chosen from
the six neighbouring vertices with the probabilities $a$, $b$ and~$c$ as in
Figure~\ref{fig:Lattice}. At the apex, that is, if $X^x_n=0$, the walk moves
to either of the vertices on row~$1$ of the lattice with probability $1/N(1)$.

If $n_L$ and~$n_R$ are both nonnegative, then it only remains to specify the
transition probabilities from the boundary vertices. However, when either
$n_L$ or~$n_R$ is negative a little more work needs to be done near the top
of the wedge, as explained in the previous subsection. In that case, we recall
the definition of row~$k_0$ from the previous subsection. For every $k=1,
\ldots,k_0-1$, we set the transition probability from each vertex in row~$k$
to each vertex in rows~$k\pm1$ equal to~$(a+c)/N(k)$. Furthermore, for every
vertex in row~$k$ that have both a left and a right neighbour, the transition
probabilities to these two neighbours are equal to~$b$. To the two vertices
on row~$k$ on the sides of~$W$ (that have only one neighbour) we assign
transition probabilities to the single neighbour and to the vertices
themselves equal to~$b$. This takes care of all the nonzero transition
probabilities near the top of the wedge.

It remains to define the transition probabilities from the boundary vertices.
As we explained in the previous subsection, we are looking for transition
probabilities that satisfy conditions~1--3 on page~\pageref{conditions}.
Below we shall give explicit expressions for these transition probabilities
from the left boundary vertices for arbitrary $n_L>0$. This is in fact
sufficient to allow us to derive the transition probabilities for all
possible wedges (i.e.~for all combinations of $n_L$ and~$n_R$), as we shall
explain first.

Observe that by left-right symmetry we can derive the transition probabilities
from the right boundary vertices for any given~$n_R$, if we know the
corresponding transition probabilities from the left boundary vertices
for $n_L=n_R$. Because the case $n_L=0$ was already covered in
Section~\ref{ssec:wedge}, it only remains to show that one can obtain the
transition probabilities from the left boundary vertices for negative~$n_L$
from those for positive~$n_L$. To do so, observe from
Figure~\ref{fig:ExtraVertices} that the left side of a wedge~$W$ with given
$n_L<0$ coincides with the right side of a (different) wedge~$W'$ with $n_R
=-n_L-1$. We propose that at the $j$th vertex on row~$k$ of~$W$ one can take
the probability of a step~$S$ equal to the probability of the step~$-S$ at
the $j$th vertex on row~$k$ of~$W'$, counted from the right side. Here we
exploit the up-down symmetry inherent in condition~1 on
page~\pageref{conditions}. It follows that it is indeed sufficient to provide
the transition probabilities from the left boundary vertices for
positive~$n_L$ only.

So let $n_L>0$ be fixed. To specify the transition probabilities we will
need some notation. We write $p_0[j,l\/]$ for the transition probability from
the $j$th vertex on a row~$k$ to the $l$th vertex on the same row. By
$p_\pm[j,l\/]$ we denote the transition probability from the $j$th vertex on
a row~$k$ to the $l$th vertex on the row $k\pm1$ (vertices on a row~$k$ are
numbered $0,1,\ldots,N(k)-1$ from left to right). If $n_R$ and~$n_L$ are both
nonnegative, then these transition probabilities are to be used at all rows
$k>0$ of the lattice, otherwise they are valid for the rows $k\geq k_0$.
Finally, we write $q_j[S]$ for the probability of a step~$S$ from the $j$th
vertex on a row~$k$. Remember from the previous subsection that there are
$N_L=n_L+1$ left boundary vertices on the rows of~$G$, so that we have to
give the transition probabilities for $j=0,1,\ldots,n_L$.

\begin{figure}
 \centering\includegraphics{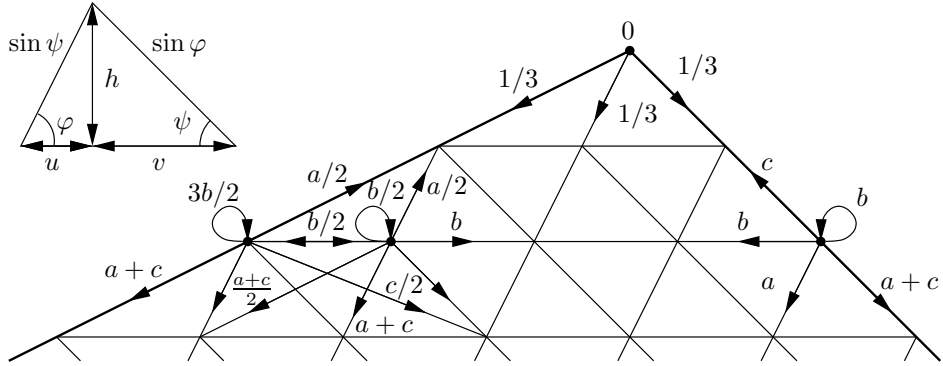}
 \caption{Transition probabilities for the reflected random walk in a wedge
  with $n_L=1$ and~$n_R=0$. The inset shows the lattice dimensions.}
 \label{fig:WedgeGeneric}
\end{figure}

First we specify the transition probabilities from a given row to the row
above. For all $j=0,1,\ldots,n_L$,
\begin{equation}
 p_-[j,0] = q_j[u+(u+v)(n_L-j)+\im h] = \frac{1}{n_L+1}a.
\end{equation}
Second, for all $j=0,1,\ldots,n_L$ the transition probabilities to the same
row are given by
\begin{eqnarray}
 && p_0[j,j+1] = q_j[u+v] = \frac{j+1}{n_L+1}b, \\
 && p_0[j,j-1] = q_j[-(u+v)] = \frac{j}{n_L+1}b, \\
 && p_0[j,j] = q_j[0] = \frac{2(n_L-j)+1}{n_L+1}b.
\end{eqnarray}
Third, specifying the transition probabilities to the row below is a bit
more complicated. For each of the boundary vertices $j=0,1,\ldots,n_L$ we have
\begin{equation}
 p_+[j,2n_L+1] = q_j[-u+(u+v)(n_L+1-j)-\im h] = \frac{c}{n_L+1}.
\end{equation}
However, for $j=0$ we have
\begin{eqnarray}
 && p_+[0,n_L] = q_j[-u-\im h] = \frac{(a+c)n_L}{n_L+1}, \\
 && p_+[0,0] = q_j[-u-(u+v)n_L-\im h] = a+c,
\end{eqnarray}
whereas for $j=1,2,\ldots,n_L-1$,
\begin{eqnarray}
 && p_+[j,j] = q_j[-u-(u+v)n_L-\im h] = \frac{(a+c)n_L}{n_L+1},
 \phantom{\qquad} \\
 && p_+[j,n_L+j] = q_j[-u-\im h] = \frac{(a+c)n_L}{n_L+1}, \\
 && p_+[j,2j] = q_j[-u-(u+v)(n_L-j)-\im h] = \frac{a+c}{n_L+1},
\end{eqnarray}
and finally, for $j=n_L$ we have
\begin{eqnarray}
 && \hskip-2em p_+[n_L,2n_L] = q_j[-u-\im h] = a+c, \\
 && \hskip-2em 
  p_+[n_L,2(n_L-n)-1] = q_j[-u-(u+v)(2n+1)-\im h] = \frac{a+c}{n_L+1}.
\end{eqnarray}
In the last equation, $n$ is an integer taking values in $\{0,1,\ldots,n_L
-1\}$. The above list specifies all the nonzero transition probabilities from
the left boundary vertices. Figure~\ref{fig:WedgeGeneric} shows an example
of the transition probabilities in the case~$n_L=1$.

We deliberately gave both the transition probabilities and the corresponding
step probabilities at the boundary vertices, to make it easy to verify that
conditions 1--3 on page~\pageref{conditions} are indeed satisfied. To verify
conditions 1 and~2, one simply has to add up the relevant transition
probabilities. To check condition~3, one may compute the real and imaginary
parts of the first step of the random walk from each of the $n_L+1$ boundary
vertices from the step probabilities (the result is given in the following
paragraph). As we explained in the previous subsection, the random walk
defined  above will converge in the scaling limit to a reflected Brownian
motion. Moreover, the walk is uniform by construction and satisfies
Lemma~\ref{lem:uniform}. Thus, the RBM will hit any horizontal cross-section
of~$W$ with the uniform distribution.

The proof of Theorem~\ref{the:main} will therefore be complete if we can
show that the angles of reflection of the RBM are $\theta_L=\alpha$
and~$\theta_R=\beta$. As in Section~\ref{sec:restricted}, the direction of
reflection with respect to the left side of~$W$ is given by the expected
value of the first step from a left boundary vertex. From the step
probabilities we compute
\begin{eqnarray}
 \Exp_z[\Re S_1] &=& \frac{1}{n_L+1} \Big[
    c(v-u)+b(u+v) \nonumber\\
 && \phantom{\frac{1}{n_L+1} \Big[}\null-(a+c)(2n_L u + n_L^2 (u+v)) \Big],\\
 \Exp_z[\Im S_1] &=& \frac{1}{n_L+1} \Big[
    -h (2c+2(a+c)n_L) \Big],
\end{eqnarray}
where $z$ is any left boundary vertex. Observing that $u/h=\cot\phi$ and
$v/h=\cot\psi$, this gives us the following expression for the reflection
angle~$\theta_L$ with respect to the left side of the wedge:
\begin{equation}
 \cot(\theta_L+\alpha)
  = \frac{\cot\phi[(a+c)(n_L^2 + 2n_L) - b + c]
      + \cot\psi[(a+c)n_L^2 - b - c]}
	 {2c+2(a+c)n_L}.
\end{equation}
Substituting Equations~(\ref{equ:a})--(\ref{equ:c}) yields
\begin{equation}
 \cot(\theta_L+\alpha) = \frac{[(\cot\phi+\cot\psi)n_L+\cot\phi]^2-1}
  {2(\cot\phi+\cot\psi)n_L+2\cot\phi}.
\end{equation}
According to Equation~(\ref{equ:alpha}), the result simplifies to
\begin{equation}
 \cot(\theta_L+\alpha) = \frac{\cot^2\alpha-1}{2\cot\alpha} = \cot(2\alpha).
\end{equation}
Thus, the reflection angle of the Brownian motion with respect to the left
side is simply $\theta_L=\alpha$. By symmetry, the angle of reflection with
respect to the right side is~$\beta$. This is consistent with the results
obtained in Section~\ref{sec:restricted} and completes the proof of
Theorem~\ref{the:main}.

\section{Properties of the RBMs}
\label{sec:properties}

This section reviews several properties shared by the reflected Brownian
motions with constant reflection angles, and sheds some light on the
RBM$_{\alpha,\beta}$ when $\alpha+\beta\geq\pi$. The properties discussed
below will be used in the following section to derive distribution functions
for the RBMs.

\subsection{Intertwining relations for the uniform walk}
\label{ssec:intertwiningdiscrete}

From the definitions of the uniform walks in Sections~\ref{sec:restricted}
and~\ref{sec:generic} the following interesting picture arises. Started with
the uniform distribution from the vertices in a given row of the lattice, the
walk can be seen as a walk from row to row on the lattice that remains uniform
on each row all the time. This can be stated more precisely in the form of an
\emph{intertwining relation}, as was noted for the case of symmetric wedges
by Dub\'edat~\cite{dubedat:2003}. Here we shall describe the generalization
of his result to generic wedges.

First let us explain what is meant by an intertwining relation. Let $(P_t:
t\geq0)$ and $(P'_t:t\geq0)$ be two Markovian semigroups with discrete or
continuous time parameter and corresponding state spaces $(S,\mathcal{S})$
and $(S',\mathcal{S'})$. Suppose that~$\Lambda$ is a Markov transition
kernel from $S'$ to~$S$. That is, $\Lambda$ is a function $\Lambda:S'\times
\mathcal{S}\to[0,1]$ such that (1)~for each fixed $x'\in S'$, $\Lambda(x',
\,\cdot\,)$ is a probability measure on~$\mathcal{S}$, and (2)~for each fixed
$A\in\mathcal{S}$, $\Lambda(\,\cdot\,,A)$ is $\mathcal{S'}$-measurable. The
two semigroups $(P_t)$ and~$(P'_t)$ are said to be intertwined by~$\Lambda$
if for all $t\geq0$ and every pair $(a',A)\in S'\times\mathcal{S}$,
\begin{equation}
 \int \Lambda(a',\dif{x})P_t(x,A) = \int P'_t(a',\dif{x'})\Lambda(x',A).
\end{equation}
In a more compact notation, the semigroups are intertwined if for all $t\geq0$
the identity $\Lambda P_t=P'_t \Lambda$ between Markov transition kernels
from $S'$ to~$S$ holds. Examples of such intertwining relations have been
studied before, see for instance~\cite{carmona:1998,rogers:1981}.

In our case, we are interested in the uniform random walk~$(X_n)$ on the
graph~$G=G_{\alpha,\beta}$ covering the wedge~$W=W_{\alpha,\beta}$, as
defined in Section~\ref{sec:generic}. Let~$(P_n)$ be its semigroup, and
let~$(P'_n)$ be the semigroup of the random walk on the set~$\N$ of
natural numbers with transition probabilities
\begin{eqnarray}
 \left\{ \begin{array}{l}
  p'(k,k) = 2b \\
  p'(k,k\pm1) = (a+c)\frac{N(k\pm1)}{N(k)}
 \end{array}\right.
 \label{equ:walkonN}
\end{eqnarray}
for $k>0$ and $p'(0,1)=1$. Observe that $p'(k,l)$ is just the conditional
probability that $X_{n+1}$ will be on row~$l$ of~$G$, given that $X_n$ is
uniformly distributed on the vertices of row~$k$.

Now consider the Markov transition kernel~$\Lambda$ from $\N$ to~$G$ such
that for each~$k\in\N$, $\Lambda(k,\,\cdot\,)$ is the uniform measure on the
vertices of row~$k$ of~$G$. Recall that for the walk~$(X_n)$, if $X_0$ is
uniformly distributed on row~$k$, the walk will stay uniform on the rows
of~$G$ afterwards (Lemma~\ref{lem:uniform}). It follows that for each
$k\in\N$, $A\subset G$ and $n\geq0$,
\begin{equation}
 \sum_{x\in G} \Lambda(k,x)P_n(x,A) = \sum_{l\in\N} P'_n(k,l)\Lambda(l,A).
 \label{equ:intertwining}
\end{equation}
Hence, we have the intertwining relation $\Lambda P_n = P'_n \Lambda$. The
reader should compare this with the statement and proof of
Lemma~\ref{lem:uniform}.

This discrete intertwining relation may be used to compute the Green function
for the uniform walk~$(X^0_n)$, killed when it reaches the row $M>0$. The
computation of this Green function gives us the expected local time
spent by the walk on a given boundary vertex of~$G$, and thereby proves
Lemma~\ref{lem:localtime}.

\begin{proof}[Lemma~\ref{lem:localtime}]
 Consider the Green function $G_M'=(I-P_M')^{-1}=\sum_{n=0}^\infty (P_M')^n$
 for the random walk on~$\N$ with transition matrix~$P_M'$ as introduced in
 Equation~(\ref{equ:walkonN}), except that now the walk is killed as soon as
 it reaches level~$M$. It is not difficult to verify that the first row of
 this Green function is given by
 \begin{equation}
  G_M'[0,k] = \left\{
  \begin{array}{ll}
    \frac{1}{(N(1)-1)(a+c)} N(k)
     \left( 1-\frac{N(k)}{N(M)} \right) & \mbox{for }k > 0; \\
	1 + \frac{a+c}{N(1)}\,g_M'[0,1]\rule{0pt}{2em} &
	 \mbox{for }k = 0.
  \end{array}
  \right.
 \end{equation}
 Now consider the Green function~$G_M$ for the uniform walk $(X^0_n)$ on the
 graph~$G$ in~$W$, killed when it reaches $\{z:\Im z = -Mh\}$. Let us denote
 the $j$th vertex on row~$k$ of~$G$ by $(j,k)$. Then, by the intertwining
 relation~(\ref{equ:intertwining}) (or by Lemma~\ref{lem:uniform}), the
 Green function~$G_M$ is related to~$G_M'$ by
 \begin{equation}
  G_M[(0,0),(j,k)] =
   \frac{1}{N(k)}G_M'[0,k] \qquad \mbox{for }j=0,1,\ldots,N(k)-1,
  \label{equ:Green}
 \end{equation}
 since the expected local time spent at vertex $(j,k)$ by the walk~$(X^0_n)$
 before it is killed at row~$M$ is the same for all $j=0,1,\ldots,N(k)-1$.
 It follows that the expected local time spent at any vertex $(j,k)$ by the
 walk before it is killed at row~$M$ is smaller than $(a+c)^{-1}$.
\end{proof}

\subsection{Intertwining relations and time reversal of the RBMs}
\label{ssec:intertwining}

In the previous subsection we considered the intertwining relation between
the uniform walk on~$G$ and a random walk on the integers. Here we will
turn our attention to the scaling limit. This time, let $(Z_t)$ be the
RBM$_{\alpha,\beta}$ in $W=W_{\alpha,\beta}$, and let $(P_t)$ be its
semigroup. Consider the Markov transition kernel~$\Lambda$ from the positive
reals~$\R_+$ to~$W$ which for each fixed $y\in\R_+$ assigns the uniform
measure to the horizontal interval $[-y\cot\alpha-\im y,y\cot\beta-\im y]$.
It is clear from the intertwining relation~(\ref{equ:intertwining}) between
the random walks, that $(P_t)$ and the semigroup of the scaling limit of the
random walk on~$\N$ will be intertwined by~$\Lambda$. It remains to identify
this scaling limit. We claim that this is a 3-dimensional Bessel process, or
in other words, it is a Brownian motion on~$\R_+$ conditioned not to hit the
origin. This generalizes Proposition~1 in Dub\'edat~\cite{dubedat:2003} to
the following statement:

\begin{theorem}
 \label{the:intertwine}
 Let $(Z_t)$, $(P_t)$ and~$\Lambda$ be as above and let $(P'_t)$ be the
 semigroup of the 3-dimensional Bessel process taking values in $\R_+$. Then
 $(P_t)$ and $(P'_t)$ are intertwined by~$\Lambda$. In particular, for all
 $y\geq0$, the process $\Big(-\Im Z^{\Lambda(y,\,\cdot\,)}_t\Big)$ is a
 3-dimensional Bessel process started from~$y$.
\end{theorem}

\begin{proof}
 As we remarked above, it is sufficient to identify the scaling limit of
 the random walk~$(X'_n)$ on~$\N$. Since the rows of the lattice have
 spacing~$h$, the proper scaling limit is obtained by considering the linear
 interpolations of the processes $h X'_{\floor{n^2 t}}/n\sigma$ where
 $\sigma^2=2(a+c)h^2$ is as before, and then taking $n\to\infty$. We may
 derive the infinitesimal generator for the limit by computing, for a
 sufficiently differentiable function~$f$ on~$\R_+$,
 \begin{eqnarray}
  \lefteqn{\hskip-2em\frac{a+c}{N(n\sigma x/h)}
   \left( f\left(x+\frac{h}{n\sigma}\right) N(n\sigma x/h + 1)
     - f\left(x-\frac{h}{n\sigma}\right) N(n\sigma x/h - 1)
   \right)} \nonumber\\
   &&\null + 2b f(x) - f(x) =
   \frac{1}{n^2}\left(\frac{1}{x}f'(x)+\frac{1}{2}f''(x) \right) + o(n^{-2}).
 \end{eqnarray}
 Here, $1/n^2$ is the time scaling, and we recognize the generator of the
 3-dimensional Bessel process (see Revuz and Yor~\cite{revuz:1991} Chapter~VI,
 \S3 and Chapter~III, Exercise~(1.15) for background on the 3-dimensional 
 Bessel process and its semigroup).
\end{proof}

We now turn our attention to the time-reversal of the RBMs. Precisely, let
$(Z_t)$ be an RBM$_{\alpha,\beta}$ in the triangle~$T=T_{\alpha,\beta}$
started from the top, and stopped when it hits~$[0,1]$. We are interested
in the time-reversal~$(\tilde{Z}_t)$ of this process. From the time-reversal
properties of the 3-dimensional Bessel process \cite[Proposition~VII(4.8)]
{revuz:1991} we know that until $(\tilde{Z}_t)$ first hits the boundary
of~$T$, it is a complex Brownian motion started with the uniform distribution
from~$[0,1]$ and conditioned not to return to~$[0,1]$. In fact, we can
describe the full process~$(\tilde{Z}_t)$ in terms of a conditioned reflected
Brownian motion. This is again a generalization of an earlier result of
Dub\'edat~\cite[Proposition~2]{dubedat:2003}:

\begin{theorem}
 \label{the:reversal}
 The time-reversal of the RBM$_{\alpha,\beta}$ in the triangle~$T=T_{\alpha,
 \beta}$, started from the top and stopped when it hits~$[0,1]$, is the
 RBM$_{\pi-\alpha,\pi-\beta}$ in~$T$ started with the uniform distribution
 from~$[0,1]$, conditioned not to return to~$[0,1]$, and killed when it hits
 the top of the triangle.
\end{theorem}

\begin{proof}
 Recall the Green function of the uniform walk~$(X^0_n)$ of
 Equation~(\ref{equ:Green}). By Nagasawa's formula (Rogers and
 Williams~\cite[III.42]{rogers:1993}) the time-reversal of this random walk
 is a Markov process with transition probabilities
 \begin{equation}
  q_M[(j,k),(l,m)] = \frac{G_M[(0,0),(l,m)]\,p[(l,m),(j,k)]}{G_M[(0,0),(j,k)]}.
 \end{equation}
 Here, the $p[(l,m),(j,k)]$ are the transition probabilities for the
 walk~$(X^0_n)$ killed at row~$M$, as specified in
 Section~\ref{ssec:walkgeneric}. We are interested in the transition
 probabilities for the reversed process in the limit $M\to\infty$. From
 the expression for the Green function it is clear that in the limit one
 gets for $k,m>0$
 \begin{equation}
  q[(j,k),(l,m)] = p[(l,m),(j,k)].
 \end{equation}
 Observe in particular that in the interior of the wedge we recover the
 transition probabilities of the original walk. Hence, the reversed walk
 converges to a Brownian motion in the interior.

 Moreover, by condition~1 on page~\pageref{conditions} it is clear that at
 every vertex on any row~$k>0$ of the lattice, the probability that the
 reversed walk will step to the row~$k+1$ is equal to the probability to
 step to the row~$k-1$. Note especially that this is not only true at
 the interior vertices but also at the boundary vertices. Therefore, the
 expected step of the random walk from any vertex on the rows $k>0$ is
 real. In particular, on the sides of~$W$ the walk must receive an average
 reflection in the real direction toward the interior of~$W$. It follows
 that the time-reversal of the reflected Brownian motion in the wedge is
 a reflected Brownian motion with reflection angles $\pi-\alpha$ and
 $\pi-\beta$ with respect to the sides of~$W$.
\end{proof}

\subsection{Reflected Brownian motions for $\alpha+\beta\geq\pi$}
\label{ssec:upsidedown}

The fact that the RBM$_{\alpha,\beta}$ is intertwined with a three-dimensional
Bessel process sheds some light on the behaviour of the reflected Brownian
motions with reflection angles satisfying $\alpha+\beta\geq\pi$. It is the
purpose of this subsection to look at these RBMs more closely. For the
duration of this subsection we will fix $\alpha,\beta\in(0,\pi)$ such that
$\alpha+\beta\geq\pi$. Then an RBM$_{\alpha,\beta}$ in the domain $T=
T_{\alpha,\beta}$ may be described by considering an RBM$_{\pi-\alpha,
\pi-\beta}$ in the wedge $W_{\pi-\alpha,\pi-\beta}$ and putting it
upside-down, as we shall explain below.

\begin{figure}
 \centering\includegraphics{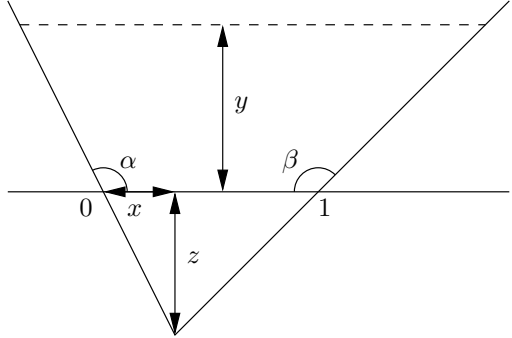}
 \caption{By putting a wedge upside-down we can shed some light on the
  RBM$_{\alpha,\beta}$ when $\alpha+\beta\geq\pi$.}
 \label{fig:UpsideDown}
\end{figure}

We set $x:=\cos\alpha\sin\beta/\sin(\alpha+\beta)$
and~$z:=-\sin\alpha\sin\beta/\sin(\alpha+\beta)$. Let $(Z_t)$ denote
RBM$_{\pi-\alpha,\pi-\beta}$ in the wedge W$_{\pi-\alpha,\pi-\beta}$, and
consider in particular the process $\Big(Z^{\Lambda(y+z,\,\cdot\,)}_t\Big)$
with $y\in\R_+$. Here, $\Lambda$ is the Markov transition kernel from~$\R_+$
to~$W_{\pi-\alpha,\pi-\beta}$ introduced in the previous subsection. Let
$f$ be the transformation $f:w\mapsto \bar{w}+x-\im z$ which puts the wedge
upside-down as illustrated in Figure~\ref{fig:UpsideDown}, and set $Y^y_t:=
f\Big(Z^{\Lambda(y+z,\,\cdot\,)}_t\Big)$. Then the process $(Y^y_t)$, stopped
when it hits the interval~$[0,1]$, is an RBM$_{\alpha,\beta}$ in $T_{\alpha,
\beta}$ started with the uniform distribution from the horizontal line
segment at altitude~$y$.

From the intertwining relation of Theorem~\ref{the:intertwine} we conclude
that the imaginary part of $(Y^y_t)$ is a Brownian motion conditioned not
to hit~$-z$. It follows (see \cite[Corollary~VI(3.4)]{revuz:1991}) that
$(Y^y_t)$ has positive probability to never reach~$[0,1]$. In fact, the
probability that $(Y^y_t)$ \emph{does} reach~$[0,1]$ is $z/(z+y)$. It is
furthermore clear from Theorem~\ref{the:intertwine} that given the event
that $(Y^y_t)$ does reach~$[0,1]$, it will arrive there with the uniform
distribution.

Now let~$\Prob_y$ denote probability with respect to the process~$(Y^y_t)$.
Then by what we said above, $(y/z+1)\Prob_y$ is a probability measure on
reflected Brownian paths in $T_{\alpha,\beta}$ started from the horizontal
line segment at altitude~$y$ that end on~$[0,1]$. Taking the limit~$y\to
\infty$ we obtain a conformally invariant probability measure on paths of
RBM$_{\alpha,\beta}$ in $T_{\alpha,\beta}$ that start ``with the uniform
distribution from infinity'' and arrive on~$[0,1]$ with the uniform
distribution. Henceforth, when we speak about RBM$_{\alpha,\beta}$ with
$\alpha+\beta\geq\pi$, we shall always assume that we restrict ourselves
to this collection of Brownian paths and the corresponding probability
measure introduced above.

\subsection{Conformal invariance and locality}
\label{ssec:invariance}

Two elementary properties shared by the RBMs are \emph{conformal invariance}
and the \emph{locality property\/}. To explain what we mean by these
properties, let us first consider an RBM$_{\alpha,\beta}$ in the wedge $W=
W_{\alpha,\beta}$ started from the point~$x\in\close{W}$. Call this process
$(Z^x_t)$, and let $v_1:=\exp((2\alpha-\pi)\im)$ and $v_2:=\exp(-2\beta\,\im)$
denote the reflection fields on the left and right sides of~$W$. Then it is
well-known (compare with Equation~(2.4) in~\cite{varadhan:1985}) that we can
uniquely write
\begin{equation}
 Z^x_t = B^x_t + v_1 Y^1_t + v_2 Y^2_t,
\end{equation}
where~$(B^x_t)$ is a complex Brownian motion started from~$x$, and $(Y^1_t)$
and~$(Y^2_t)$ are real-valued continuous increasing processes adapted
to~$(B^x_t)$ such that $Y^1_0=Y^2_0=0$. Moreover, $Y^1_t$ increases only
when~$Z^x_t$ is on the left side of~$W$ and $Y^2_t$ increases only when~$Z^x_t$
is on the right side of~$W$ ($Y^1_t$ and~$Y^2_t$ are essentially the local
times of~$(Z^x_t)$ on the two sides of~$W$).

Now let~$g$ be a conformal transformation from~$W$ onto a domain~$D$ with
smooth boundary. Consider the sum
\begin{equation}
 g(Z^x_t)-g(Z^x_0) = \sum_{j=0}^{N-1}
  \Big[ g(Z^x_{(j+1)t/N}) - g(Z^x_{jt/N}) \Big].
\end{equation}
To compute the sum we use Taylor's theorem to expand each term. The
computation is very similar to the one in Section~\ref{ssec:scalinglimit}.
In particular, on the boundary we may keep only the first-order terms. Then
letting $N\to\infty$ and using the fact that the real and imaginary parts
of~$g$ are harmonic, just as in the proof of It\^o's formula (see
e.g.~Sections 4.2 and~4.3 in Gardiner~\cite{gardiner:1983} for a nice
discussion) one obtains
\begin{equation}
 g(Z^x_t)-g(Z^x_0) = \int_0^t g'(Z^x_s)\,\dif{B^x_s}
  + v_1 \int_0^t g'(Z^x_s)\,\dif{Y^1_s} + v_2 \int_0^t g'(Z^x_s)\,\dif{Y^2_s}.
 \label{equ:conformalRBM}
\end{equation}

\begin{figure}
 \centering\includegraphics{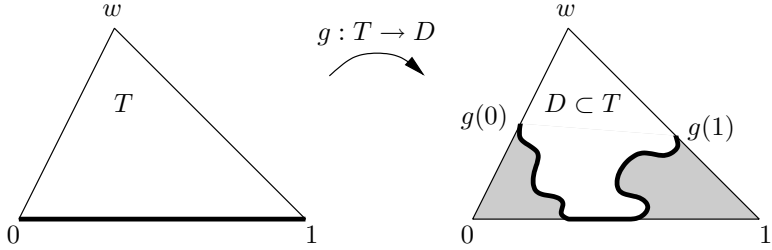}
 \caption{The locality property says that the RBM$_{\alpha,\beta}$ in~$T$
  started from~$w$ and stopped when it exits from the subset~$D$ behaves
  just like an RBM$_{\alpha,\beta}$ in~$D$ started from~$w$.}
 \label{fig:Locality}
\end{figure}

The first integral in~(\ref{equ:conformalRBM}) is the usual expression for
the conformal image of Brownian motion. Making the usual time-change
$u(t):=\int_0^t|g'(Z^x_s)|^2\,\dif{s}$ (see Revuz and Yor~\cite[Theorem~V(2.5)]
{revuz:1991}) and denoting its inverse by~$t(u)$, we conclude from
Equation~(\ref{equ:conformalRBM}) that the process $\tilde{Z}_{u} :=
g(Z^x_{t(u)})$ is a reflected Brownian motion in~$D$ with reflection vector
fields $v_1\,g'(g^{-1}(\,\cdot\,))$ and $v_2\,g'(g^{-1}(\,\cdot\,))$ on the
two ``sides'' of~$D$ (i.e.~the images of the two sides of~$W$). Note in
particular that because~$g$ is an angle-preserving transformation, the process
$(\tilde{Z}_{u})$ is also reflected at the angles $\alpha$ and~$\beta$ with
respect to the boundary of~$D$. This shows that the RBM$_{\alpha,\beta}$ is
conformally invariant.

We may use the same reasoning to explain what we mean by the locality
property of the RBM$_{\alpha,\beta}$. We change the setup to one that will be
more useful later. Indeed, we now let $(Z_t)$ be an RBM$_{\alpha,\beta}$ in
the triangle $T=T_{\alpha,\beta}$ started from the top $w=w_{\alpha,\beta}$,
and set $\tau:=\inf\{t\geq 0: Z_t\in[0,1]\}$. Furthermore, we take~$g$ to be
a conformal map of~$T$ onto an open subset~$D$ of~$T$ that fixes~$w$, and
such that the left and right sides of~$T$ are mapped onto subsets of
themselves. Finally, we set $\sigma:=\inf\{t\geq0: Z_t\in\close{T}\setminus
D\}$, the exit time of the RBM$_{\alpha,\beta}$ from the subset~$D$. See
Figure~\ref{fig:Locality} for an illustration.

From the calculation above we conclude that up to the stopping time~$\tau$,
the process $g(Z_t)$ is just a time-changed RBM$_{\alpha,\beta}$ in the subset
$D$ of the triangle~$T$. In particular, modulo a time-change the laws of
$(g(Z_t):t\leq\tau)$ and $(Z_t:t\leq\sigma)$ are the same. This is what is
called the locality property. For more background and for consequences of
the locality property we refer to the SLE literature~\cite{lsw:I,werner:2002}.

\section{Distribution functions}
\label{sec:distributions}

In this section we compute several distribution functions associated with
the family RBM$_{\alpha,\beta}$ of reflected Brownian motions. We fix
$\alpha$ and~$\beta$ in $(0,\pi)$ for the duration of the section, with no
further restrictions on $\alpha,\beta$ (recall that when $\alpha+\beta\geq
\pi$, we assume that we work with the probability measure of
Section~\ref{ssec:upsidedown}). Furthermore, we fix two angles $\lambda,\mu
\in(0,\pi)$ such that $\lambda+\mu<\pi$. These two angles define the domain
$T=T_{\lambda,\mu}$ in this section. To be somewhat more general, we will
derive distribution functions for the RBM$_{\alpha,\beta}$ in the triangle
$T=T_{\lambda,\mu}$ rather than in the domain $T_{\alpha,\beta}$.

\subsection{Notations}
\label{ssec:hullnotations}

\begin{figure}
 \centering\includegraphics{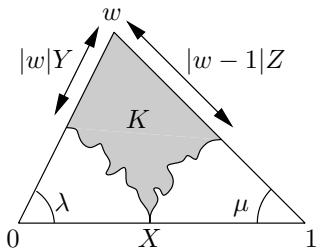}
 \caption{Definition of the hull~$K$ of the reflected Brownian motion with
  parameters $\alpha$ and~$\beta$ in the triangle~$T$, and of the random
  variables $X$, $Y$ and~$Z$.}
 \label{fig:Hull}
\end{figure}

First we introduce some notations. We shall call the set of points
disconnected from $\{0,1\}$ in the triangle~$T$ by the path of the RBM up
to the time when it hits~$[0,1]$ the \emph{hull}~$K$ of the process. This
hull has exactly one point in common with the real line, which we denote
by~$X$. We shall denote by $|w|Y$ the distance of the lowest point of the
hull on the left side to the top~$w=w_{\lambda,\mu}$, and by~$|w-1|Z$ the
distance of the lowest point of the hull on the right side to the top~$w$.
Thus, all three random variables $X$, $Y$ and~$Z$ take on values in the
range~$[0,1]$. See Figure~\ref{fig:Hull} for an illustration of the
definitions.

Below we shall compute the marginal and joint distributions of the
variables $X$, $Y$ and~$Z$. We shall see that these can be expressed in
terms of the conformal transformations of the upper half of the complex
plane onto the triangles $T_{\gamma,\delta}$. Thus it is useful to
review some properties of these transformations first. We simplify the
notation somewhat by writing $\gamma'$ as an abbreviation for $\gamma/\pi$
whenever~$\gamma$ denotes an angle.

\begin{figure}
 \centering\includegraphics{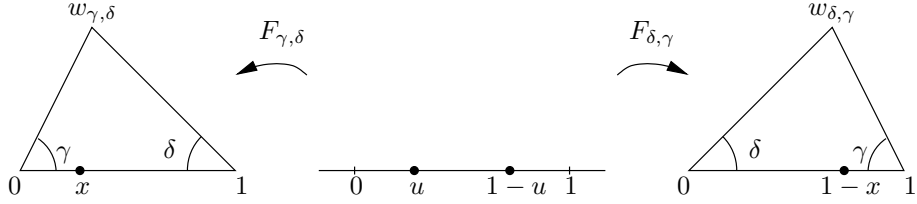}
 \caption{Transformations of the upper half-plane onto triangles.}
 \label{fig:Mapping}
\end{figure}

By the Schwarz-Christoffel formula of complex analysis
(see~\cite[Section~6.2.4]{ahlfors:1966} or~\cite[Section XI.3]{gamelin:2000}),
the unique conformal transformation of the upper half-plane~$\H$
onto~$T_{\gamma,\delta}$ that fixes $0$ and~$1$ and maps $\infty$
to~$w_{\gamma,\delta}$ is given by
\begin{equation}
 F_{\gamma,\delta}(z) = \frac{\int_0^z t^{\gamma'-1}(1-t)^{\delta'-1}\dif{t}}
   {\int_0^1 t^{\gamma'-1}(1-t)^{\delta'-1}\dif{t}}
  = \frac{B_z(\gamma',\delta')}{B(\gamma',\delta')},
 \label{equ:trianglemapping}
\end{equation}
where $B(\gamma',\delta')=\Gamma(\gamma')\Gamma(\delta')/
\Gamma(\gamma'+\delta')$ is the beta-function, and $B_z(\gamma',\delta')$
the incomplete beta-function (see e.g.~\cite{gamma} for background on these
special functions). Symmetry considerations show that the transformations
satisfy
\begin{equation}
 F_{\gamma,\delta}(u) = 1-F_{\delta,\gamma}(1-u) \quad\mbox{and}\quad
 F^{-1}_{\gamma,\delta}(x) = 1-F^{-1}_{\delta,\gamma}(1-x).
 \label{equ:symmetry}
\end{equation}
See Figure~\ref{fig:Mapping} for an illustration.

A different kind of distribution that we can compute is the conditional
probability that the last side of the triangle visited by the RBM, given
that it lands at $X=x$, is the right side. As we shall see, this
distribution can also be expressed in term of triangle mappings. In fact,
it turns out that there is a remarkable resemblance between this conditional
probability and the marginal distribution functions of the variables $X$,
$Y$ and~$Z$.

\subsection{Characteristics of the hull}
\label{ssec:hulldistributions}

Here we derive the (joint) distribution functions of $X$, $Y$ and~$Z$,
which are characteristics of the hull generated by the RBM. Remember that
we are considering an RBM$_{\alpha,\beta}$ in the triangle~$T=T_{\lambda,\mu}$
started from $w=w_{\lambda,\mu}$. For convenience, let us also introduce the
angle $\nu:=\pi-\lambda-\mu$. Then our main conclusion may be formulated as
follows:

\begin{proposition}
 Let $a=a(y)$ and~$b=b(z)$ be the points on the left and right sides of~$T$
 at distances $|w|y$ and~$|w-1|z$ from~$w$, respectively. Then the joint
 distribution of $X$, $Y$ and~$Z$ is given by
 \begin{eqnarray}
  \lefteqn{ \Prob[X\leq x,Y\leq y,Z\leq z]\phantom{\int} = \null} \\
  && F_{\alpha,\beta}\left(
      \frac{F^{-1}_{\lambda,\mu}(x)-F^{-1}_{\lambda,\mu}(a)}
 	      {F^{-1}_{\lambda,\mu}(b)-F^{-1}_{\lambda,\mu}(a)}
      \right)
 	- F_{\alpha,\beta}\left(
      \frac{-F^{-1}_{\lambda,\mu}(a)}
 	      {F^{-1}_{\lambda,\mu}(b)-F^{-1}_{\lambda,\mu}(a)}
      \right) \nonumber
 \end{eqnarray}
 where the images of $a$ and~$b$ under the map $F^{-1}_{\lambda,\mu}$ can be
 expressed in terms of $y$ and~$z$ as
 \begin{equation}
  F^{-1}_{\lambda,\mu}(a)
   = 1 - \frac{1}{F^{-1}_{\nu,\lambda}(y)};
 \end{equation}
 \begin{equation}
  F^{-1}_{\lambda,\mu}(b)
   = \frac{1}{1-F^{-1}_{\mu,\nu}(1-z)}
   = \frac{1}{F^{-1}_{\nu,\mu}(z)}.
 \end{equation}
 \label{pro:jointdist}
\end{proposition}

\noindent
Note that in the last equation we used the symmetry
property~(\ref{equ:symmetry}).

\begin{figure}
 \centering\includegraphics{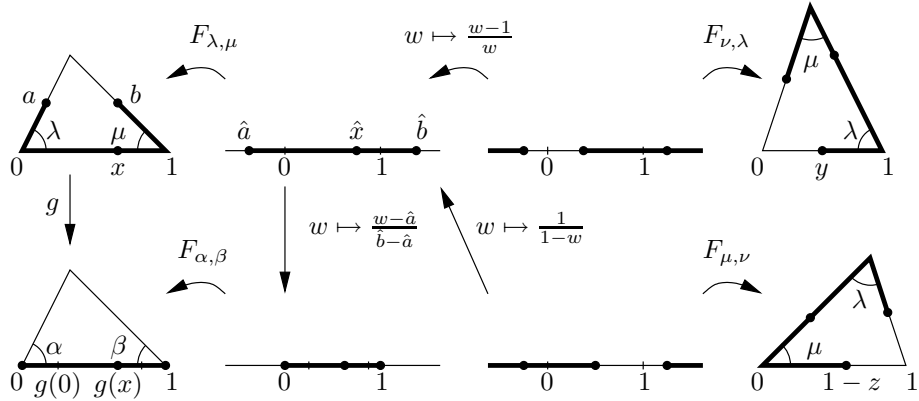}
 \caption{This illustration shows schematically how one obtains an explicit
  form for the map~$g$ in terms of the variables $y$ and~$z$. The notations
  $\ahat$, $\bhat$ and~$\xhat$ in the figure are short for $F^{-1}_{\lambda,
  \mu}(a)$, $F^{-1}_{\lambda,\mu}(b)$ and~$F^{-1}_{\lambda,\mu}(x)$.}
 \label{fig:JointBig}
\end{figure}

\begin{proof}
 The idea of the computation of $\Prob[X\leq x,Y\leq y,Z\leq z]$ is
 illustrated in Figure~\ref{fig:JointSmall}. Consider an RBM$_{\alpha,\beta}$
 in the triangle~$T$ started from the top~$w$, and stopped as soon as it hits
 the counter-clockwise arc from $a$ to~$b$ on the boundary (the thick line
 in the figure). Then the probability $\Prob[X\leq x,Y\leq y,Z\leq z]$ is
 just the probability that this process is stopped in the interval~$(0,x)$.

\begin{figure}
 \centering\includegraphics{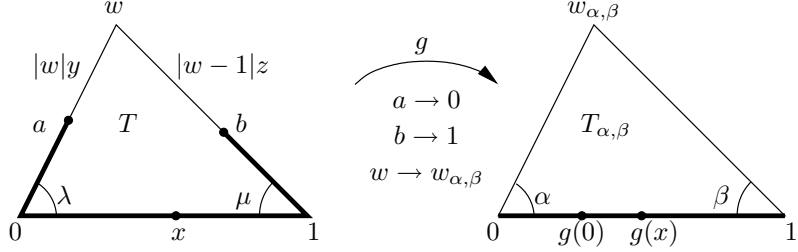}
 \caption{This figure illustrates how the joint distribution function of the
  random variables $X$, $Y$ and~$Z$ can be computed. As explained in the text,
  the joint probability $\Prob[X\leq x,Y\leq y,Z\leq z]$ is just $g(x)-g(0)$.}
 \label{fig:JointSmall}
\end{figure}

 We now use conformal invariance and locality. Let $g=g_{a(y),b(z)}$ be the
 conformal map of~$T$ onto~$T_{\alpha,\beta}$ that sends $a$ to~$0$, $b$
 to~$1$ and $w$ to~$w_{\alpha,\beta}$, as illustrated in
 Figure~\ref{fig:JointSmall}. Then the probability that we are trying to
 compute is exactly the probability that an RBM$_{\alpha,\beta}$
 in~$T_{\alpha,\beta}$ started from~$w_{\alpha,\beta}$ and stopped when it
 hits~$[0,1]$, is stopped in the interval~$(g(0),g(x))$. But since the exit
 distribution of the RBM is uniform in~$T_{\alpha,\beta}$, this probability
 is simply $g(x)-g(0)$. Thus,
 \begin{equation}
  \Prob[X\leq x,Y\leq y,Z\leq z] = g(x)-g(0).
 \end{equation}
 Our next task is to find an explicit expression for this joint probability
 by deriving the explicit form of the map~$g=g_{a(y),b(z)}$. The explicit
 form of~$g$ is obtained by suitably combining conformal self-maps of the
 upper half-plane with triangle mappings. How this is done exactly is
 described in Figure~\ref{fig:JointBig}. The expression for the joint
 distribution follows.
\end{proof}

By sending one or or two of the three variables $x$, $y$ and~$z$ to~$1$, and
using the symmetry property~(\ref{equ:symmetry}), one may derive the following
corollaries of Proposition~\ref{pro:jointdist}:

\begin{corollary}
 We have the following joint distributions:
 \begin{eqnarray}
  \lefteqn{\Prob[X\leq x, Z\leq z] = F_{\alpha,\beta}\left(
   F^{-1}_{\lambda,\mu}(x) F^{-1}_{\nu,\mu}(z) \right);} \\
  \lefteqn{\Prob[X\leq x, Y\leq y] =
   F_{\beta,\alpha}\left( F^{-1}_{\nu,\lambda}(y) \right) -
   F_{\beta,\alpha}\left(
    F^{-1}_{\mu,\lambda}(1-x) F^{-1}_{\nu,\lambda}(y) \right);} \\
  \lefteqn{\Prob[Y\leq y,Z\leq z] =
   F_{\alpha,\beta}\left( \frac{F^{-1}_{\nu,\mu}(z)}
    {F^{-1}_{\nu,\mu}(z) + F^{-1}_{\nu,\lambda}(y)
     - F^{-1}_{\nu,\mu}(z)F^{-1}_{\nu,\lambda}(y)}
    \right)} \\
   &&\null + F_{\beta,\alpha}\left( \frac{F^{-1}_{\nu,\lambda}(y)}
    {F^{-1}_{\nu,\mu}(z) + F^{-1}_{\nu,\lambda}(y)
     - F^{-1}_{\nu,\mu}(z)F^{-1}_{\nu,\lambda}(y)}
    \right) - 1.\hskip6em \nonumber
 \end{eqnarray}
 \label{cor:jointdist}
\end{corollary}

\begin{corollary}
 The marginal distributions of $X$, $Y$ and~$Z$ are given by
 \begin{eqnarray}
  \Prob[X\leq x] &=& F_{\alpha,\beta}\left( F^{-1}_{\lambda,\mu}(x) \right); \\
  \Prob[Y\leq y] &=& F_{\beta,\alpha}\left( F^{-1}_{\nu,\lambda}(y) \right);\\
  \Prob[Z\leq z] &=& F_{\alpha,\beta}\left( F^{-1}_{\nu,\mu}(z) \right).
 \end{eqnarray}
 \label{cor:marginaldist}
\end{corollary}

Observe that the marginal distributions of the variables $X$, $Y$ and~$Z$
take on particularly simple forms. These marginal distribution functions have
a nice geometric interpretation. For instance, $\Prob[Y\leq y]$ is just the
image of~$y$ under the transformation that maps the triangle $T_{\nu,\lambda}$
onto~$T_{\beta,\alpha}$, fixes $0$ and~$1$ and takes $w_{\nu,\lambda}$
onto~$w_{\beta,\alpha}$. Similar observations hold for $\Prob[X\leq x]$
and $\Prob[Z\leq z]$.

These observations lead to some intriguing conclusions. First of all, we
conclude that for an RBM$_{\pi/3,\pi/3}$ in an equilateral triangle (so
that $\alpha=\beta=\lambda=\mu=\nu=\pi/3$) all three variables $X$, $Y$
and~$Z$ are uniform. This is not so surprising when we realize that the
hull in this case is the same as that of the exploration process of
critical percolation, as we noted in the introduction. Indeed, in the case
of percolation $X$ and~$Y$ can be interpreted as the endpoints of the
highest crossing of a given colour between the sides $(0,w)$ and~$(0,1)$
of the triangle. Thus by symmetry, if $X$ is uniform, then so are $Y$
and~$Z$.

\begin{figure}
 \centering\includegraphics{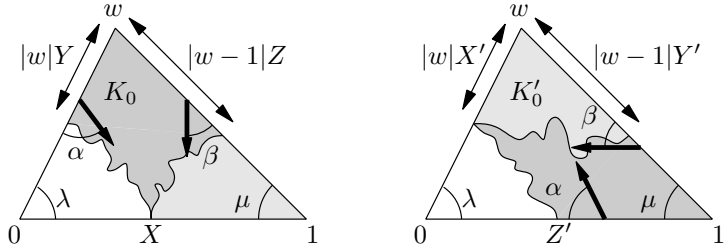}
 \caption{The unions $K_0$ and~$K'_0$ of the shaded sets on the left and
  right are generated by different RBMs, but have the same law.}
 \label{fig:Intriguing}
\end{figure}

In other triangles, similar but more intricate connections exist. For
example, let~$(Z_t)$ be an RBM$_{\alpha,\beta}$ in~$T=T_{\lambda,\mu}$
started from~$w=w_{\lambda,\mu}$ and stopped when it hits~$[0,1]$, as
before. Compare this process with an RBM~$(Z'_t)$ in~$T$ started from~$1$,
stopped on~$[0,w]$ and reflected on $(0,1)$ at an angle~$\alpha$ and on
$(w,1)$ at an angle~$\beta$ with respect to the boundary. Here it is assumed
that small angles denote reflection away from~$1$. To this second RBM we can
associate normalized random variables $X'$, $Y'$ and~$Z'$, measuring the
distances of the exit point on $[0,w]$ to~$w$ and of the ``lowest points''
of the hull on $[w,1]$ and~$[0,1]$ to $w$ and~$0$, respectively. See
Figure~\ref{fig:Intriguing}. It follows from Corollary~\ref{cor:marginaldist}
that $X$ and~$Z'$ have the same distribution (and so do $Y$ and~$X'$).

This result has an interesting interpretation in terms of the hulls generated
by the two processes, as we shall now describe. We write $\Omega$ for the
collection of closed, connected subsets~$C$ of $\close{T}$ such that the
right side of~$T$ is in~$C$ and $T\setminus C$ is connected. We further
define~$\mathcal{Q}$ as the collection of compact $A\subset\close{T}$ such
that $A=\close{A\cap T}$, and $\close{T}\setminus A$ is simply connected and
contains the right side of~$T$. We then endow~$\Omega$ with the $\sigma$-field
generated by the events $\{K\in\Omega:K\cap A=\emptyset\}$ for all
$A\in\mathcal{Q}$. This setup is similar to the one in Lawler, Schramm
and Werner~\cite[Sections 2 and~3]{lsw:2003}. In particular, a probability
measure~$\Prob$ on~$\Omega$ is characterized by the values of $\Prob[K\cap A
=\emptyset]$ for $A\in\mathcal{Q}$, see~\cite[Lemma~3.2]{lsw:2003}.

\begin{theorem}
 Consider the processes $(Z_t)$ and~$(Z'_t)$ stopped on $[0,1]$ and $[0,w]$
 as described above. Let $K_0$ and~$K'_0$ denote the sets of points in
 $\close{T}$ that are disconnected from~$0$ by the paths of $(Z_t)$ and
 $(Z'_t)$, respectively. Then the laws of $K_0$ and~$K'_0$ on the
 space~$\Omega$ are the same.
\end{theorem}

\begin{proof}
 Let~$\Prob$ be the law of~$K_0$, and let~$A\in\mathcal{Q}$. Denote by $a$
 and~$b$ the points of $A\cap\partial T$ closest to~$w$ and~$1$, respectively.
 Let $g:T\to T\setminus A$ be the conformal transformation that fixes $1$
 and~$w$, and maps~$0$ onto~$a$ if $\Im a>0$, and onto~$0$ otherwise. Then,
 by conformal invariance of the RBM$_{\alpha,\beta}$ in~$T$, $\Prob[K_0\cap
 A=\emptyset]=\Prob[X\in(g^{-1}(b),1)]$ (compare this with our discussion of
 the locality property in Section~\ref{ssec:invariance}). Likewise, the
 law~$\Prob'$ of~$K'_0$ satisfies $\Prob'[K'_0\cap A=\emptyset]=\Prob'
 [Z'\in(g^{-1}(b),1)]$. Since $X$ and~$Z'$ have the same distribution by
 Corollary~\ref{cor:marginaldist}, the theorem follows.
\end{proof}

\subsection{Last-visit distribution}
\label{ssec:lastvisit}

As we discussed in Section~\ref{ssec:intertwining}, the imaginary part of
the RBM$_{\alpha,\beta}$ in $T_{\alpha,\beta}$ stopped when it hits~$[0,1]$
is a three-dimensional Bessel process, and so is the imaginary part of its
time-reversal. In particular, the time-reversed process, considered up to
the first time when it hits the left or right side of the triangle, is
a complex Brownian motion started with uniform distribution from~$[0,1]$
and conditioned not to return to the real line. In other words, up to its
first contact with the left or right side, this process is a  Brownian
excursion of the upper half-plane, started with uniform distribution
from~$[0,1]$ (for background on Brownian excursions, see \cite{lsw:I}
and~\cite{lawler:1999b}). This fact allows us to derive the following
result:

\begin{proposition}
 Let $(Z_t)$ be an RBM$_{\alpha,\beta}$ in the triangle $T=T_{\lambda,\mu}$
 started from~$w=w_{\lambda,\mu}$. Let~$\tau:=\inf\{t\geq0:Z_t\in[0,1]\}$,
 and let~$\sigma$ be the last time before~$\tau$ when~$(Z_t)$ visited the
 boundary of~$T$. Let $E$ denote the event that $Z_\sigma$ is on the right
 side of the triangle. Then
 \[ \Prob[E\mid Z_\tau=x]
  = F_{\pi-\alpha,\pi-\beta}\left( F^{-1}_{\lambda,\mu}(x) \right).
 \]
\end{proposition}

\noindent
It is shown in Dub\'edat~\cite{dubedat:2003} how one derives this result in
the special case where $\alpha=\beta=\lambda=\mu=\pi/3$. Here, for the sake
of completeness, we repeat the computation for the general case.

\begin{figure}
 \centering\includegraphics{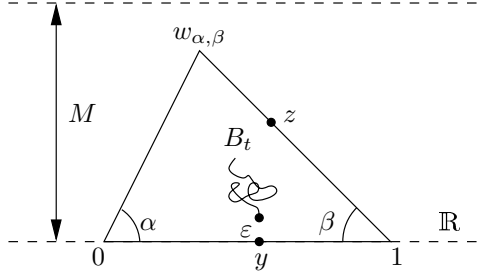}
 \caption{A complex Brownian motion~$B_t$ started from $y+\im\epsilon$ in
  the triangle~$T_{\alpha,\beta}$. This process conditioned to exit the strip
  $\{z : 0<\Im z<M\}$ through the top boundary and stopped when it hits the
  boundary of the triangle~$T_{\alpha,\beta}$ corresponds in the limit
  $M\to\infty$, $\epsilon\downarrow0$ to the time-reversal of the
  RBM$_{\alpha,\beta}$ in the triangle, as explained in the text.}
 \label{fig:LastVisit}
\end{figure}

\begin{proof}
 We want to use the fact that the time-reversal of the RBM$_{\alpha,\beta}$
 in the triangle $T_{\alpha,\beta}$ starts out as a Brownian excursion. So
 we first map $T_{\lambda,\mu}$ onto~$T_{\alpha,\beta}$ by the transformation
 $F_{\alpha,\beta}\circ F^{-1}_{\lambda,\mu}$. This maps the point~$x$ onto
 $y:=F_{\alpha,\beta}\left(F^{-1}_{\lambda,\mu}(x)\right)$.

 Next, let $(Y_t:t\geq0)$ be a Brownian excursion of the upper half-plane,
 and let $S$ be the first time when~$(Y_t)$ visits either the left or right
 side of~$T_{\alpha,\beta}$. Then
 \begin{equation}
  \Prob[E\mid Z_\tau=x] = \lim_{\epsilon\downarrow0}\int_1^{w_{\alpha,\beta}}
   \Prob_{y+\im\epsilon}[Y_S\in\dif{z}],
  \label{equ:lastvisit1}
 \end{equation}
 where the integrals runs over the right side of the triangle~$T_{\alpha,
 \beta}$, and~$\Prob_z$ denotes probability with respect to the Brownian
 excursion started from~$z$.

 Now let $(B_t:t\geq0)$ be a complex Brownian motion, let~$U$ be the first
 time when~$(B_t)$ visits either the left or right side of~$T_{\alpha,\beta}$,
 and let~$U_M$ be the first time when~$(B_t)$ exits the strip $\{z:
 0<\Im z<M\}$. Suppose that $\Prob'_z$ denotes probability with respect to
 the Brownian motion started from~$z$. Then, using the strong Markov property
 of Brownian motion (see Figure~\ref{fig:LastVisit} for a sketch), we have
 \begin{eqnarray}
  \Prob_{y+\im\epsilon}[Y_S\in\dif{z}] &=&
   \lim_{M\to\infty} \Prob'_{y+\im\epsilon}[B_U\in\dif{z} \mid \Im B_{U_M} = M]
   \nonumber\\
  &=& \Prob'_{y+\im\epsilon}[B_U\in\dif{z}] \lim_{M\to\infty}
   \frac{\Prob'_z[\Im B_{U_M} = M]}{\Prob'_{y+\im\epsilon}[\Im B_{U_M} = M]}
   \nonumber\\
  &=& \Prob'_{y+\im\epsilon}[B_U\in\dif{z}]\frac{\Im z}{\epsilon},
  \label{equ:lastvisit2}
 \end{eqnarray}
 where in the last step we have used~\cite[Proposition~II(3.8)]{revuz:1991}.
 Combining Equations (\ref{equ:lastvisit1}) and~(\ref{equ:lastvisit2}), we
 see that we have to compute the limit of $\Prob'_{y+\im\epsilon}[B_U\in
 \dif{z}]\Im z/\epsilon$ as $\epsilon\to0$. This computation can be done by
 using the conformal invariance of Brownian motion.

\begin{figure}
 \centering\includegraphics{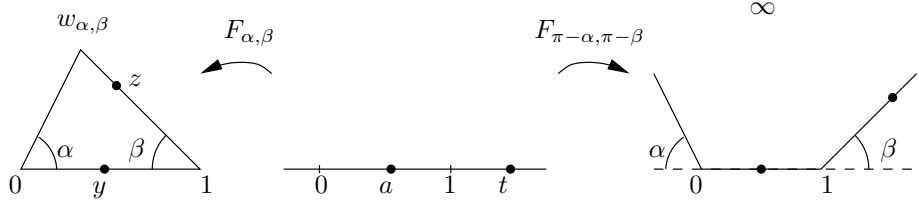}
 \caption{Conformal transformations between the upper half-plane and the
  domains $T_{\alpha,\beta}$ and~$T_{\pi-\alpha,\pi-\beta}$.}
 \label{fig:LastVisitMap}
\end{figure}

 Remember that the probability that a complex Brownian motion started from
 $a+\im b$ leaves the upper half-plane through~$(-\infty,x)$, is given by
 the harmonic measure~$\omega(x)$ of~$(-\infty,x)$ at the point $a+\im b$. It
 is straightforward to verify that
 \begin{equation}
  \omega(x) = \int_{-\infty}^x \frac{b}{\pi}\frac{\dif{t}}{b^2+(t-a)^2}
   = \frac{1}{2} + \frac{1}{\pi} \arctan\frac{x-a}{b}.
 \end{equation}
 Thus, mapping the triangle~$T_{\alpha,\beta}$ conformally to the upper
 half-plane by the transformation~$F^{-1}_{\alpha,\beta}$ as in
 Figure~\ref{fig:LastVisitMap}, we can write
 \begin{equation}
  \lim_{\epsilon\downarrow0}\Prob'_{y+\im\epsilon}[B_U\in\dif{z}]
   \frac{\Im z}{\epsilon}
   = \frac{\Im z}{\pi F_{\alpha,\beta}'(a)}
     \frac{(F^{-1}_{\alpha,\beta})'(z)\dif{z}}
 	 {(F^{-1}_{\alpha,\beta}(z)-a)^2}
 \end{equation}
 where $a=F^{-1}_{\alpha,\beta}(y)=F^{-1}_{\lambda,\mu}(x)$. Therefore,
 using~(\ref{equ:trianglemapping}),
 \begin{eqnarray}
  \lefteqn{\Prob[E \mid Z_\tau=x]} \\
   &&\null
    = \frac{1}{\pi F_{\alpha,\beta}'(a)} \int_1^\infty \frac{\dif{t}}{(t-a)^2}
    \Im F_{\alpha,\beta}(t) \nonumber\\
   &&\null
    = \frac{a^{1-\alpha'}(1-a)^{1-\beta'}}{\pi} \int_1^\infty
    \frac{\dif{t}}{(t-a)^2}
    \Im \int_1^t u^{\alpha'-1}(1-u)^{\beta'-1}\dif{u} \nonumber\\
   &&\null
    = \frac{\sin\beta}{\pi} a^{1-\alpha'}(1-a)^{1-\beta'} \int_1^\infty
    \dif{u}\,u^{\alpha'-1}(u-1)^{\beta'-1} \int_u^\infty \frac{\dif{t}}{(t-a)^2}
    \nonumber\\
   &&\null
    = \frac{\sin\beta}{\pi} a^{1-\alpha'}(1-a)^{1-\beta'} \int_1^\infty
    u^{\alpha'-1}(u-1)^{\beta'-1}(u-a)^{-1} \dif{u} \nonumber\\
   &&\null
    = \frac{\sin\beta}{\pi} a^{1-\alpha'}(1-a)^{1-\beta'} \int_0^1
    t^{1-\alpha'-\beta'}(1-t)^{\beta'-1}(1-at)^{-1} \dif{t}, \nonumber
 \end{eqnarray}
 where in the last step we have made the substitution~$t=1/u$.

 Using equations 15.3.1 and~15.2.5 for the hypergeometric function
 in~\cite{hypergeometric} and the formulas 6.1.15 and~6.1.17 for the gamma
 function from~\cite{gamma}, we finally derive
 \begin{eqnarray}
  \lefteqn{\Prob[E\mid Z_\tau=x]} \\
   && \null = \frac{\sin\beta}{\pi}
    \frac{\Gamma(2{-}\alpha'{-}\beta')\Gamma(\beta')}{\Gamma(2{-}\alpha')}
    a^{1-\alpha'}(1-a)^{1-\beta'}
    \Hypergeom(1,2{-}\alpha'{-}\beta';2{-}\alpha';a)
    \nonumber\\
   &&\null = \frac{\Gamma(2-\alpha'-\beta')}{\Gamma(1-\alpha')\Gamma(1-\beta')}
    \int_0^a t^{-\alpha'} (1-t)^{-\beta'} \dif{t} \nonumber\\
   &&\null = F_{\pi-\alpha,\pi-\beta}(a)
   = F_{\pi-\alpha,\pi-\beta}\left( F^{-1}_{\lambda,\mu}(x)
    \right).\phantom{\int_0^0} \nonumber
 \end{eqnarray}
 This is the desired result.
\end{proof}

In words, we have considered the conditional probability that the last side
visited by an RBM$_{\alpha,\beta}$ in~$T_{\lambda,\mu}$ started
from~$w_{\lambda,\mu}$ is the right side, given that the exit point $X$
equals~$x\in(0,1)$. This conditional probability is exactly given by the
image of~$x$ under the transformation that maps $T_{\lambda,\mu}$
onto~$T_{\pi-\alpha,\pi-\beta}$, fixing $0$ and~$1$ and sending~$w_{\lambda,
\mu}$ onto~$w_{\pi-\alpha,\pi-\beta}$.

\paragraph{Acknowledgements.}
The author wishes to thank Julien Dub\'edat for a useful discussion on the
subject of this paper. Thanks are also due to Remco van der Hofstad for
his comments and aid in preparing the manuscript. This research was supported
financially by the Dutch research foundation FOM (Fundamenteel Onderzoek der
Materie).

\end{document}